\documentclass[12pt]{article}
\usepackage{amssymb,amsfonts,amsmath, psfrag,eepic,colordvi,graphicx,epsfig,ytableau}
\usepackage{amssymb,latexsym,graphics,array}
\usepackage{MnSymbol}
\usepackage[enableskew]{youngtab}
\usepackage{float,enumerate,enumitem}
\usepackage{tikz,ifthen}
\usepackage{pgfplots}
\usetikzlibrary{math}
\topskip  
\newcommand{\rmnum}[1]{\romannumeral #1}

\numberwithin{equation}{section}
\newtheorem{theorem}{Theorem}[section]

\newtheorem{corollary}[theorem]{Corollary}

\newtheorem{lemma}[theorem]{Lemma}

\newtheorem{observation}[theorem]{Observation}

\begin{document}
	\parskip 6pt
	
	\pagenumbering{arabic}
	\def\sof{\hfill\rule{2mm}{2mm}}
	\def\ls{\leq}
	\def\gs{\geq}
	\def\SS{\mathcal S}
	\def\qq{{\bold q}}
	\def\MM{\mathcal M}
	\def\TT{\mathcal T}
	\def\EE{\mathcal E}
	\def\lsp{\mbox{lsp}}
	\def\rsp{\mbox{rsp}}
	\def\pf{\noindent {\it Proof.} }
	\def\mp{\mbox{pyramid}}
	\def\mb{\mbox{block}}
	\def\mc{\mbox{cross}}
	\def\qed{\hfill \rule{4pt}{7pt}}
	\def\block{\hfill \rule{5pt}{5pt}}
	\def\lr#1{\multicolumn{1}{|@{\hspace{.6ex}}c@{\hspace{.6ex}}|}{\raisebox{-.3ex}{$#1$}}}

	\begin{center}
		{\Large \bf  Equidistribution of  set-valued statistics  on standard Young tableaux and transversals}
	\end{center}
	
	\begin{center}
		{\small  Robin D.P. Zhou$^a$,  Sherry H.F. Yan$^{b,*}$\footnote{$^*$Corresponding author.}  \footnote{{\em E-mail address:} hfy@zjnu.cn. }}
		
		$^a$College of Mathematics Physics and Information\\
		Shaoxing University\\
		Shaoxing 312000, P.R. China
		
		$^{b}$Department of Mathematics,
		Zhejiang Normal University\\
		Jinhua 321004, P.R. China		
	\end{center}
	
	\noindent {\bf Abstract.}  As a natural generalization  of permutations, transversals of Young diagrams play an important role in the study of pattern avoiding permutations.  
	Let $\mathcal{T}_{\lambda}(\tau)$ and  $\mathcal{ST}_{\lambda}(\tau)$ denote   the set of $\tau$-avoiding transversals and $\tau$-avoiding  symmetric transversals
	of a Young diagram $\lambda$, respectively.  
	In this paper, we are mainly concerned with the distribution of  the peak set and the valley set  on standard Young tableaux and pattern avoiding transversals. In particular, by introducing Knuth transformations on standard Young tableaux, we prove that  the peak set and the valley set are equidistributed on   the  standard Young tableaux of  shape  $\lambda/\mu$ for any skew diagram $\lambda/\mu$.  
	The  equidistribution enables us to show that   the peak set is equidistributed over $\mathcal{T}_{\lambda}(12\cdots k\tau)$ (resp.  $\mathcal{ST}_{\lambda}(12\cdots k \tau)$)  and  $\mathcal{T}_{\lambda}(k\cdots 21\tau) $   (resp. $\mathcal{ST}_{\lambda}(k\cdots 21\tau)$) for any Young diagram $\lambda$ and any permutation $\tau$ of $\{k+1, k+2, \ldots, k+m\}$ with $k,m\geq 1$. 
	Our results  are refinements of the result  of  Backelin-West-Xin which states that $|\mathcal{T}_{\lambda}(12\cdots k\tau)|=|\mathcal{T}_{\lambda}(k\cdots 21  \tau)|$ and the result of Bousquet-M\'elou  and Steingr\'imsson  which states that $|\mathcal{ST}_{\lambda}(12\cdots k \tau)|=|\mathcal{ST}_{\lambda}(k\cdots 21  \tau)|$. 
	As    applications,  we are able to  
	\begin{itemize}
		\item confirm a recent conjecture    posed by  Yan-Wang-Zhou which asserts   that  the peak set is equidistributed over  $  12\cdots k  \tau$-avoiding involutions  and  $  k\cdots 21  \tau$-avoiding involutions;
		\item prove that 
		alternating involutions avoiding the pattern $  12\cdots k  \tau$ are equinumerous with  alternating involutions avoiding the pattern $  k\cdots 21 \tau$,  paralleling the  result of Backelin-West-Xin for  permutations,  the result of Bousquet-M\'elou  and Steingr\'imsson   for involutions,   and the result of Yan for alternating permutations.  
	\end{itemize}



	
	\section{Introduction}
	As a natural generalization  of permutations, transversals of Young diagrams play an important role in the study of pattern avoiding permutations and various interesting results
	on them  have been  obtained  in the literature (see e.g. \cite{BW, BWX, Bousquet, Dukes, Jaggard, Yan2023, Zhou}).  This paper is devoted to the investigation of  the distribution of  the peak set and the valley set  on standard Young tableaux and pattern avoiding transversals.

	Let us first review some necessary definitions before we  state our original motivation and  main results.
	Let $\mathcal{S}_n$ denote the set of permutations of [n]=\{1,2,\ldots,n\}, 
	which we will always write as words $\pi = \pi_1\pi_2\cdots \pi_n$. 
	An index $i$ $(2 \leq i \leq n-1)$ is said to be a {\em peak} (resp. {\em valley}) of $\pi$ if 
	$\pi_{i-1} < \pi_i > \pi_{i+1}$ (resp. $\pi_{i-1} > \pi_i < \pi_{i+1}$).  
	Let $\mathrm{Peak}(\pi)$ and $\mathrm{Val}(\pi)$ denote the set of peaks and the set of valleys of $\pi$, respectively.  For example, let $\pi = 5\,6\,1\,9\,4\,3\,7\,2\,8$. 
	Then we have  $\mathrm{Peak}(\pi) = \{2,4,7\}$ and $\mathrm{Val}(\pi) = \{3,6,8\}$.
	There  are many research articles devoted to the combinatorics of 
	peaks on permutations  (see e.g. \cite{ Aguiar,Billera,Petersen,Schocker,Stanley,Stembridge,Strehl,Warren}).

	Let $\pi = \pi_1\pi_2\cdots \pi_n\in \mathcal{S}_n$.
	A permutation $\pi$ is said to be {\em alternating} if $\pi_1 < \pi_2 > \pi_3 < \pi_4 > \cdots$.
	A permutation $\pi$ is said to be an {\em involution} if $\pi= \pi^{-1}$, where $\pi^{-1}$ denotes the inverse permutation of $\pi$.
	For example, the permutation $\pi = 7\,9\,5\,6\,3\,4\,1\,8\,2$ is an alternating involution.
	Let $\mathcal{A}_n$ (resp.  $\mathcal{I}_n$ and  $\mathcal{AI}_n$) denote the set of alternating permutations (resp.   involutions and alternating involutions) of length $n$.

	Given a permutation $\pi \in \mathcal{S}_n$ and a permutation $\sigma \in \mathcal{S}_k$,
	an occurrence of $\sigma$ in $\pi$ is a subsequence $\pi_{i_1}\pi_{i_2}\cdots \pi_{i_k}$
	of $\pi$ that is order isomorphic to $\sigma$.
	We say $\pi$   contains  the pattern $\sigma$ if $\pi$ contains an occurrence of $\sigma$.
	Otherwise, we say $\pi$ avoids the pattern $\sigma$  and
	$\pi$ is $\sigma$-avoiding. For instance,  the permutation $5\,7\,2\,4\,3\,6\,1\,8$ is $1243$-avoiding while it contains a pattern $1234$
	corresponding to the subsequence $2468$.
	Let $\mathcal{S}_n(\sigma)$ denote the set of  $\sigma$-avoiding permutations of length $n$.
	We will keep this notation also when $\mathcal{S}_n$ is replaced by
	other subsets of permutations such as $\mathcal{A}_n, \mathcal{I}_n, \mathcal{AI}_n$,  etc.

	Pattern avoiding permutations were introduced by Knuth \cite{Knuth} in 1970 and first
	systematically studied by Simion-Schmidt \cite{Simion}.
	The theory of pattern avoidance has been extensively exploited  over particular subsets of permutations.
	For example, various results have been obtained for
	pattern avoiding  alternating permutations
	(see e.g. \cite{Bona,Chen,Lewis2009,Lewis2011,Lewis2012,Mansour2003,Stanley,Yan2012,Yan2013})  and pattern avoiding involutions (see e.g. \cite{Barnabei2011,Bona2016,Bousquet,Dukes,Guibert,Jaggard}).

	Barnabei-Bonetti-Castronuovo-Silimbani \cite{Barnabei} initiated  the enumeration of pattern avoiding alternating involutions  \cite{Barnabei}.
	They enumerated and characterized some classes of alternating   involutions avoiding a
	single pattern of length $4$. 
	Furthermore, Barnabei-Bonetti-Castronuovo-Silimbani  \cite{Barnabei} posed several
	conjectures concerning pattern avoiding alternating   involutions, which have been confirmed  by Yan-Wang-Zhou \cite{Yan2023} and Zhou-Yan \cite{Zhou}.

	Given an integer $n$, a {\em partition} of $n$ is a sequence of
	nonnegative integers $\lambda=(\lambda_1,\lambda_2,\ldots,\lambda_k)$
	such that $n=\lambda_1+\lambda_2+\cdots+\lambda_k$ and $\lambda_1 \geq \lambda_2 \geq \cdots \geq \lambda_k >0$.
	We denote $\lambda \vdash n$ or $|\lambda| = n$.
	Here $k$ is called the length of $\lambda$, denoted by $\ell(\lambda)$.
	A {\em Young diagram} of shape $\lambda$ is defined to be a left-justified array of $n$ squares with $\lambda_1$ squares
	in the first row, $\lambda_2$ squares in the second row and so on.
	The {\em conjugate} of a partition $\lambda$, denoted by $\lambda^{T}$, is the partition whose Young diagram is the reflection along the main diagonal of $\lambda$'s Young diagram, and $\lambda$ is said to be {\em self-conjugate} if $\lambda=\lambda^{T}$.
	In this paper, a square $(i,j)$ in a Young diagram $\lambda$ is referred to the square
	in the $i$-th column and $j$-th row of $\lambda$,
	where we number the rows from top to bottom, and the columns from left to right. In what follows, we always
	treat  a partition $\lambda$ and its Young diagram as 
	identical.  
	
	Let $\mu$ and $\lambda$ be two Young diagrams with  $\mu \subseteq \lambda$ (i.e. $\mu_i \leq \lambda_i$ for all $i$) and $|\lambda| - |\mu| = n$.
	The {\em skew diagram} $\lambda/\mu$ is  defined to be the diagram obtained from the Young diagram
	$\lambda$ by removing the Young diagram $\mu$ at the top-left corner. In this context, 
	$n$ is said to be the size of $\lambda/\mu$.  
	In what follows, we treat $\lambda/\emptyset$ and 
	$\lambda$ as identical.
	A {\em standard Young tableau} (SYT) $T$ of (skew) shape $\lambda/ \mu$ is a filling of the skew diagram $\lambda/\mu$
	with the numbers
	$1,2,\ldots,n$ such that the entries are increasing alone each rows and
	each columns.
	Figure \ref{fig:SYT} (left) illustrates an SYT of shape $(4,3,2)$ and 
	Figure \ref{fig:SYT} (right) illustrates an SYT of
	shape $(5,4,2) / (2,1)$.
	Let  $\mathrm{SYT}(\lambda /\mu)$
	denote the set of SYT's of  
	shape $\lambda /\mu$. 
	Let $f^{\lambda/\mu}$ denote the cardinality of $\mathrm{SYT}(\lambda/\mu)$.
	Note that $f^{\lambda/\mu}$ can be given by the famous hook-length formula.

	\begin{figure} [h]
		\centering
		{$\begin{array}[b]{*{4}c}\cline{1-4}
				\lr{1}&\lr{4}&\lr{6}&\lr{9}\\\cline{1-4}
				\lr{2}&\lr{5}&\lr{8}\\\cline{1-3}
				\lr{3}&\lr{7}\\\cline{1-2}
			\end{array}$}
		\hskip 2.5cm
		$\begin{array}[b]{*{5}c}\cline{3-5}
			&&\lr{1}&\lr{2}&\lr{6}\\\cline{2-5}
			&\lr{3}&\lr{4}&\lr{7}\\\cline{1-4}
			\lr{5}&\lr{8}\\\cline{1-2}
		\end{array}$
		\caption{Examples of standard Young tableaux.}\label{fig:SYT}
	\end{figure}
	
	For an SYT $T$ of  shape $\lambda /\mu$ with $n$ entries,
	an index $i$ ($1\leq i \leq n-1$) is said to be
	a {\em descent}  of $T$ if $i + 1$ appears in a lower row of $T$ than $i$, 
	otherwise, $i$ is said to be  an  {\em ascent} of $T$.
	An index $i$  $(2\leq i\leq n-1)$ is said to be a {\em peak} of $T$  if $i-1$ is an ascent and $i$ is a descent, whereas an index $i$  $(2\leq i\leq n-1)$ is said to be a {\em valley} of  $T$  if $i-1$ is a descent and $i$ is an ascent. 
	Let $\mathrm{Peak}(T)$ and $\mathrm{Val}(T)$ denote the set of peaks and the set of valleys of the SYT $T$, respectively. 
	For example, let $T$ be the SYT as shown in Figure \ref{fig:SYT} (left), 
	we have $\mathrm{Peak}(T) = \{4,6\}$ and $\mathrm{Val}(T) = \{3,5,7\}$.
	
	The Robinson-Schensted-Knuth algorithm (RSK algorithm) \cite{StanleyVol2} sets up a bijection between symmetric group $\mathcal{S}_n$ and pairs $(P,Q)$ of SYT's of the same shape $\lambda \vdash n$.
	We denote this correspondence by $\pi \xrightarrow{\mathrm{RSK}} (P,Q)$, where $\pi \in \mathcal{S}_n$.
	We call $P$ the {\em insertion tableau} and $Q$ the {\em recording tableau} of $\pi$.
	The {\em shape} of $\pi$, denoted by $\mathrm{sh}(\pi)$, is defined to be the shape of $P$ (or $Q$).
	By the insertion rule of the RSK algorithm, 
	it can be easily seen that $\pi$ has the same peak set and the same valley set as its recording tableau $Q$.

	For a set $A$, write
	$t^A = \prod\limits_{i\in A}t_i $.
	Set $t^A = 1$ when $A = \emptyset$.
	Our first main result is concerned with the equidistribution of the set-valued statistics $\mathrm{Peak}$  and $\mathrm{Val}$ on the  SYT's of the same shape.

	\begin{theorem}\label{thm:skewSYT}
		Let $n\geq 1$. 
		Then for  any  skew diagram  $\lambda /\mu$  of size $n$,  
		the set-valued statistics $\mathrm{Peak}$  and $\mathrm{Val}$  are equidistributed over
		$\mathrm{SYT}(\lambda /\mu)$,  that is,   
		$$\sum_{T \in \mathrm{SYT}(\lambda /\mu) }t^{\mathrm{Peak}(T)} = \sum_{T \in \mathrm{SYT}(\lambda /\mu) }t^{\mathrm{Val}(T)}.$$
	\end{theorem}

	It turns out that Theorem   \ref{thm:skewSYT} can be used 
	to study the distribution of the set-valued statistics on pattern avoiding transversals. This is actually the motivation and original intention behind our writing of this article. The equiditribution of permutation statistics on pattern avoiding permutations and transversals   has been investigated in recent literature, see \cite{Bloom1, Bloom2, Bloom3,  Chan, dd, Yan2015, Yan2023, Zhou} and references therein.   
	
\begin{figure}[H]
		\begin{center}
			\begin{tikzpicture}[font =\small , scale = 0.3, line width = 0.7pt]
				\foreach \i / \j in {1/9,2/9,3/9,4/9,5/7,6/7,7/7,8/4,9/4}
				{
					\foreach \k in {1,...,\j} \draw (\i,-\k)rectangle(\i + 1,-\k -1);
				}
				\foreach \i / \j in {1/2,2/1,3/8,4/9,5/7,6/6,7/5,8/3,9/4} \filldraw[black](\i+0.5,-\j-0.5)circle(5pt);
			\end{tikzpicture}
		\end{center}
		\caption{A transversal of the Young diagram $\lambda=(9,9,9,9,7,7,7,4,4)$.}\label{fig:transversal}
	\end{figure}
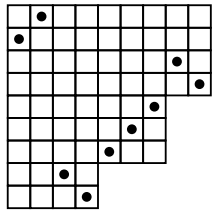

	A {\em transversal} $T$ of shape $\lambda$ is a $01$-filling of the squares of the  Young diagram $\lambda$ such that each row and column contains exactly one $1$.
	A permutation $\pi = \pi_1\pi_2\cdots\pi_n$ can be regarded as a transversal of the $n$ by $n$ square diagram,
	in which the square $(i,\pi_i)$ is filled with a $1$ for all $1\leq i\leq n$
	and all the other squares  are filled with $0$'s.
	The transversal corresponding to the permutation $\pi$ is also called the
	{\em permutation matrix} of $\pi$.
	Let $\mathcal{T}_{\lambda}$ denote the set of transversals of the Young diagram $\lambda$
	and let $\mathcal{T}_n$ denote the set of transversals of all Young diagrams with $n$ columns.
	A transversal $T$ of shape $\lambda$ will be written  as $T = t_1t_2\cdots t_n$
	if the square $(i,t_i)$ in $T$ contains a $1$ for $1\leq i \leq n$ and other squares of $T$ 
	contains $0$'s. 
	For example, Figure \ref{fig:transversal} illustrates the transversal $T = 2\,1\,8\,9\,7\,6\,5\,3\,4$ of shape $\lambda = (9,9,9,9,7,7,7,4,4)$
	in which each $\bullet$ represents a $1$ and each empty square represents a 
	$0$.
	
	A transversal $T = t_1t_2\cdots t_n$ of shape $\lambda$ is said to be {\em symmetric} if $\lambda$ is self-conjugate and the corresponding permutation $t_1t_2\cdots t_n$ is an involution. More precisely,  in a   symmetric  transversal $T$ of a self-conjugate Young diagram $\lambda$,  the square $(i,j)$ contains a $1$ if and only if the square $(j,i)$ contains a $1$. 
	The transversal $T$ as shown in 
	Figure \ref{fig:transversal} gives an example of  symmetric transversals.
	Let $\mathcal{ST}_{\lambda}$ denote the set of symmetric transversals of the self-conjugate Young diagram $\lambda$
	and let $\mathcal{ST}_n$ denote the set of symmetric transversals of all self-conjugate Young diagrams with $n$ columns.

	For a Young diagram $\lambda$, let $c_j(\lambda)$ denote the length of 
	the $j$-th column of $\lambda$. 
	Let $T$ be a transversal of shape $\lambda$.
	An index $i$ ($2 \leq i \leq  n-1$) with $c_{i-1}(\lambda) = c_{i}(\lambda) = c_{i+1}(\lambda)$ is said to be a {\em peak} (resp. valley) of 
	$T$ if $t_{i-1}< t_i > t_{i+1}$ (resp. $t_{i-1}> t_i < t_{i+1}$).
	Let $\mathrm{Peak}(T)$ and $\mathrm{Val}(T)$ denote the set of peaks 
	and the set of valleys of $T$, respectively. It is easily seen that when restricted to permutations, the peaks (resp. valleys) of transversals coincide with the peaks (resp. valleys) of permutations.

	The notion of pattern avoidance on permutations can be extended to transversals.
	Let $T$ be a transversal of shape $\lambda$.
	We say $T$ contains the
	pattern $\sigma$ if there exist $R=\{r_1<r_2<\cdots<r_k\}$ and $C=\{c_1< c_2< \cdots< c_k\}$,
	such that each of the squares $(c_i,r_j)$ falls within the board of $\lambda$ and
	the matrix  restricted on rows $R$ and columns $C$ is the permutation matrix of $\sigma$.
	Otherwise, we say $T$ avoids $\sigma$ and $T$ is $\sigma$-avoiding.
	Let $\mathcal{T}_n(\sigma)$ denote the set of $\sigma$-avoiding transversals
	in $\mathcal{T}_n$.
	We will keep this notation also when $\mathcal{T}_n$ is replaced by
	other subsets of $\mathcal{T}_n$.
	For example, the transversal as shown in Figure \ref{fig:transversal}
	is $4321$-avoiding while it contains the pattern $2134$.
	
	Given two permutations $\pi = \pi_1\pi_2\cdots \pi_n \in \mathcal{S}_n$
	and $\sigma = \sigma_1\sigma_2\cdots \sigma_m \in \mathcal{S}_m$,
	the {\em direct sum} of $\pi$ and $\sigma$, denoted by $\pi \oplus \sigma$, is the permutation $\pi_1\pi_2\cdots \pi_n (\sigma_1+n)(\sigma_2+n)\cdots (\sigma_m+n)$.
	Let $\epsilon$ denote the empty permutation.
	Throughout the paper, we treat $\pi \oplus \epsilon$ and 
	$\pi$ as identical.
	Let $I_k = 12\cdots k$ and $J_k = k\cdots21$. 
	Our  second main result is concerned with the following  equidistribution of set-valued statistics on pattern avoiding transversals. 
	
	\begin{theorem}\label{thm:transversal}
		For any Young diagram $\lambda$ and any positive integer $k$, we have 
		$$\sum_{T \in \mathcal{T}_{\lambda}(I_k) }t^{\mathrm{Peak}(T)}
		= \sum_{T \in \mathcal{T}_{\lambda}(I_k) }t^{\mathrm{Val}(T)}
		= \sum_{T \in \mathcal{T}_{\lambda}(J_k) }t^{\mathrm{Val}(T)}
		= \sum_{T \in \mathcal{T}_{\lambda}(J_k) }t^{\mathrm{Peak}(T)}.$$
	\end{theorem}
	
	\begin{theorem}\label{thm:transversal-general}
		Let $k\geq 1$. For any Young diagram $\lambda$ and any permutation $\tau$, the set-valued statistic $\mathrm{Peak}$  is equidistributed over $\mathcal{T}_{\lambda}(I_k\oplus \tau)$  and $\mathcal{T}_{\lambda}(J_k\oplus \tau)$, that is, 
		$$\sum_{T \in \mathcal{T}_{\lambda}(I_k\oplus \tau) }t^{\mathrm{Peak}(T)}
		= \sum_{T \in \mathcal{T}_{\lambda}(J_k\oplus \tau) }t^{\mathrm{Peak}(T)}.$$
	\end{theorem}

	Setting $t_i=1$ for all $i\geq 1$   in Theorem \ref{thm:transversal-general} recovers the result of Backelin-West-Xin \cite{BWX} which states that $|\mathcal{T}_{\lambda}(I_k\oplus \tau)|=|\mathcal{T}_{\lambda}(J_k\oplus \tau)|$.  
	Hence, Theorem \ref{thm:transversal-general}  can be viewed as a refinement of the result of Backelin-West-Xin.  
	We remark that some examples show that 
	the bijection between $\mathcal{T}_{\lambda}(I_k\oplus \tau)$ and $\mathcal{T}_{\lambda}(J_k\oplus \tau)$ established by Backelin-West-Xin  does not preserve the peak set.

	Recall that when we set $\lambda$ to be the $n$ by $n$ square diagram,
	the transversal $T\in \mathcal{T}_{\lambda}$ becomes a 
	permutation $\pi$ in $\mathcal{S}_n$ with  $\mathrm{Peak}(\pi) = \mathrm{Peak}(T)$ and  $\mathrm{Val}(\pi) = \mathrm{Val}(T)$.  
	Then the following result follows directly from Theorem  \ref{thm:transversal-general}. 
	\begin{corollary} 
		Let $n, k\geq 1$. For any permutation $\tau$, the peak set is equidistributed over   $\mathcal{S}_n(I_k\oplus \tau )$ and $\mathcal{S}_n(J_k\oplus \tau )$. 
	\end{corollary}

	For symmetric transversals, we have the following analogues of Theorems
	\ref{thm:transversal} and \ref{thm:transversal-general}.
	
	\begin{theorem}\label{thm:symm-transversal}
		Let $k\geq 1$. For any self-conjugate Young diagram $\lambda$, we have 
		$$\sum_{T \in \mathcal{ST}_{\lambda}(I_k) }t^{\mathrm{Peak}(T)}
		= \sum_{T \in \mathcal{ST}_{\lambda}(I_k) }t^{\mathrm{Val}(T)}
		= \sum_{T \in \mathcal{ST}_{\lambda}(J_k) }t^{\mathrm{Val}(T)}
		= \sum_{T \in \mathcal{ST}_{\lambda}(J_k) }t^{\mathrm{Peak}(T)}.$$
	\end{theorem}
	
	\begin{theorem}\label{thm:symm-transversal-general}
		Let $k\geq 1$. For any self-conjugate Young diagram $\lambda$ and any pattern $\tau$, the set-valued statistic $\mathrm{Peak}$  is equidistributed over $\mathcal{ST}_{\lambda}(I_k\oplus \tau)  $ and  $\mathcal{ST}_{\lambda}(J_k\oplus \tau)  $, that is, 
		$$\sum_{T \in \mathcal{ST}_{\lambda}(I_k\oplus \tau) }t^{\mathrm{Peak}(T)}
		= \sum_{T \in \mathcal{ST}_{\lambda}(J_k\oplus \tau) }t^{\mathrm{Peak}(T)}.$$
	\end{theorem}
	
	Note that the  case $k=3$ of Theorem  \ref{thm:symm-transversal-general}  has been verified by  Yan-Wang-Zhou \cite{Yan2023}   by establishing a peak set preserving bijection between  $\mathcal{ST}_{\lambda}(J_3\oplus \tau)$ and $\mathcal{ST}_{\lambda}(I_3\oplus \tau)$.   

Setting   $t_i=1$ for all $i\geq 1$   in Theorem \ref{thm:symm-transversal-general} recovers the result of 
	Bousquet-M\'elou  and Steingr\'imsson \cite{Bousquet} which states that $|\mathcal{ST}_{\lambda}(I_k\oplus \tau)|=|\mathcal{ST}_{\lambda}(J_k\oplus \tau)|$.  Therefore, Theorem \ref{thm:symm-transversal-general} can be viewed as a refinement of the result of 
	Bousquet-M\'elou  and Steingr\'imsson. 
	We remark that some examples show that 
	the bijection between $\mathcal{ST}_{\lambda}(I_k\oplus \tau)$ and $\mathcal{ST}_{\lambda}(J_k\oplus \tau)$ established by Bousquet-M\'elou  and Steingr\'imsson  does not preserve the peak set. 
	

	Recall that when we set $\lambda$ to be the $n$ by $n$ square diagram,
	the transversal $T\in \mathcal{ST}_{\lambda}$ becomes an
	involution $\pi$ in $\mathcal{I}_n$ with  $\mathrm{Peak}(\pi) = \mathrm{Peak}(T)$ and  $\mathrm{Val}(\pi) = \mathrm{Val}(T)$. Hence, the following result follows directly from Theorem \ref{thm:symm-transversal-general}, confirming a recent conjecture posed by Yan-Wang-Zhou \cite{Yan2023}.
	
	\begin{corollary}(\cite{Yan2023}, Conjecture 4.1)
		Let $n,k\geq 1$. For any permutation $\tau$,  the set-valued statistic $\mathrm{Peak}$ is equidistributed over  $\mathcal{I}_n(I_k\oplus \tau)$  and  $\mathcal{I}_n (J_k\oplus \tau)$, that is, 
		$$\sum_{\pi \in \mathcal{I}_n(I_k\oplus \tau) }t^{\mathrm{Peak}(\pi)}
		= \sum_{\pi\in \mathcal{I}_n (J_k\oplus \tau) }t^{\mathrm{Peak}(\pi)}.$$
	\end{corollary}

	Recently,  Yan-Wang-Zhou \cite{Yan2023} proved that $|\mathcal{AI}_n(I_3 \oplus \tau)| =|\mathcal{AI}_n(J_3 \oplus \tau)|$  for any  nonempty permutation $\tau$ as conjectured by Barnabei-Bonetti-Castronuovo-Silimbani \cite{Barnabei}. 
	 In this paper, we shall obtain the following  extension of   the result of Yan-Wang-Zhou to general $k$. 
	\begin{theorem}\label{thm:AI}
		Let $n,k\geq 1$. For any nonempty permutation $\tau$, 
		we have $$|\mathcal{AI}_n(I_k \oplus \tau)| =|\mathcal{AI}_n(J_k \oplus \tau)|.$$ 
	\end{theorem}
	Note that Backelin-West-Xin \cite{BWX}  proved that  $|\mathcal{S}_n(I_k \oplus \tau)| =|\mathcal{S}_n(J_k \oplus \tau)|$, which has been extended to involutions by  Bousquet-M\'elou  and Steingr\'imsson  \cite{Bousquet} and to alternating permutations by Yan \cite{Yan2013}. Hence, Theorem \ref{thm:AI} can be  viewed as  a  parallel  work of  the above results.

	The rest of the  paper is organized as follows.
	Section 2 is devoted to the investigation of  the equidistribution  of the set-valued  statistics  Peak   and  Val 
	on standard Young tableaux of a skew shape.
	By introducing Knuth transformations on standard Young tableaux, we prove Theorem
	\ref{thm:skewSYT}.
	In Section \ref{sec:3}, relying on Theorem \ref{thm:skewSYT}, we shall investigate the distribution of Peak and Val on   transversals, thereby proving Theorems \ref{thm:transversal}, \ref{thm:transversal-general},
	\ref{thm:symm-transversal}, \ref{thm:symm-transversal-general} and  \ref{thm:AI}.

	\section{Peaks and valleys on standard Young tableaux}\label{sec:2}

	This section is  devoted to the proof of 
	Theorem \ref{thm:skewSYT}.
	To this end, we shall associate each SYT  with a permutation and show that 
	all the permutations corresponding to the SYT's of a given skew shape can be 	divided into several Knuth equivalence classes.

		Given an SYT $T$ of shape $\lambda/\mu$ with $n$ entries,  we associate
	$T$ with a word $y = y_1y_2\cdots y_n = \alpha(T)$ of length $n$, where
	$y_i$ is the row index of the square of $T$ containing the number $i$.
	The word $y$ is known as {\em Yamanouchi word}. In this context, we say that $y$ is a Yamanouchi word {\em compatible} to 
	$\lambda/\mu$.
	
	For a word $w$, we denote   $|w|_i$ to be the  the number of occurrences of $i$ in $w$.
	Let $\lambda=(\lambda_1, \lambda_2, \ldots, \lambda_k)$ and $\mu=(\mu_1, \mu_2, \ldots, \mu_{\ell})$ be two partitions satisfying that $k \geq \ell$ and $\lambda_i \geq \mu_i$ for all $i$ with the assumption $\mu_i = 0$ when $i > \ell$.
	In fact, a word $y=y_1y_2\cdots y_n$ is a  Yamanouchi word compatible to $\lambda/\mu$ if and only if 
	\begin{itemize}
		\item  $y_j \in [k]$ for all $1\leq j\leq n$;
		\item $|y|_i=\lambda_i-\mu_i$ for all $1\leq i\leq k$;
		\item for all $1\leq i\leq k-1$ and $1\leq j\leq n$, we have $|y^{(j)}|_{i}+\mu_{i}\geq |y^{(j)}|_{i+1}+\mu_{i+1} $,  where $y^{(j)}=y_1y_2\cdots y_j$.
	\end{itemize}
	Denote by $\mathcal{Y}(\lambda/\mu)$ the set of Yamanouchi words compatible to the skew diagram $\lambda/\mu$. 
	
	By the above analysis, it is straightforward to  recover the corresponding SYT of 
	shape $\lambda/\mu$ from its the  Yamanouchi word $y$
	by letting the $i$-th row of $T$
	contain the indices of the letters of $y$ which are equal to $i$.
	For example, the Yamanouchi word of the SYT $T$
	in Figure \ref{fig:SYT} (left) is given by  $123121321$ and the Yamanouchi word of the SYT $S$ in Figure \ref{fig:SYT} (right) is given by $11223123$. 
	
	For a Yamanouchi word $y   = y_1y_2\cdots y_n$ of an SYT $T$, 
	an index $i$ ($2\leq i \leq n-1$)  is said to be a {\em peak} of $y$ if $y_{i-1} < y_i \geq y_{i+1}$, and  $i$ is said to be a {\em valley} of $y$ if $y_{i-1} \geq  y_i < y_{i+1}$.
	Let $\mathrm{Peak}(y)$ and $\mathrm{Val}(y)$ denote the set of peaks and the set of valleys of the Yamanouchi word $y$, respectively.
 The following result follows directly from  the definition of the Yamanouchi word of an SYT. 
	
	\begin{lemma}\label{lem:Yamanouchi-word}
		For any skew diagram $\lambda/\mu$, the map $\alpha$ induces a bijection between $\mathrm{SYT}(\lambda/\mu)$ and $\mathcal{Y}(\lambda/\mu)$ such that for any $T\in \mathrm{SYT}(\lambda/\mu)$, we have 
         $\mathrm{Peak}(T) = \mathrm{Val}(\alpha(T))$ and 
		$\mathrm{Val}(T) = \mathrm{Peak}(\alpha(T))$.
	\end{lemma}
	
	
	Let $\lambda /\mu$ be a skew diagram of size $n$ such that $\lambda_i - \mu_i = a_i$ for $1\leq i \leq \ell(\lambda)$ with the convention  $\mu_i = 0$ for $i > \ell(\mu)$.
	Given an SYT $T$ of skew shape $\lambda / \mu$,  let $y = y_1y_2\cdots y_n$ 
	be its corresponding Yamanouchi word.
	We associate $T$ with a permutation  
	$\pi = \pi_1\pi_2\cdots \pi_n = \beta(T)$ of length $n$, where 
	$\pi$ is obtained from $y$ be converting the $a_i$ $i$'s in $y$ to the numbers 
	$a_1+a_2+\cdots +a_{i-1}+1$, $a_1+a_2+\cdots +a_{i-1}+2,\ldots$, $a_1+a_2+\cdots +a_{i-1}+a_i$ from left to right in decreasing order.
	We call the permutation $\pi$ the {\em Yamanouchi permutation}.
	It is easily seen that $\beta$ gives an injective map from 
	$\mathrm{SYT}(\lambda /\mu)$ to $\mathcal{S}_n$.
	On the other hand, given a Yamanouchi permutation $\pi$, we can recover the Yamanouchi word $y\in \mathcal{Y}(\lambda/\mu)$ from $\pi$ by converting $a_1+a_2+\cdots +a_{i-1}+1$, $a_1+a_2+\cdots +a_{i-1}+2,\ldots$, $a_1+a_2+\cdots +a_{i-1}+a_i$ in $\pi$ to $i$,
	thereby recovering  the SYT $T$ of shape $\lambda /\mu$ from 
	the Yamanouchi word $y$.
	For example, the Yamanouchi permutation of the SYT $T$ 
	in Figure \ref{fig:SYT} (left) is given by $479362851$ and the Yamanouchi permutation of the SYT $S$ in Figure \ref{fig:SYT} (right) is given by $32658147$.
	One can easily check that $\mathrm{Peak}(y) = \mathrm{Peak}(\pi)$ and 
	$\mathrm{Val}(y) = \mathrm{Val}(\pi)$.
	Combining Lemma \ref{lem:Yamanouchi-word}, we immediatey obtain the following  result. 
	
	\begin{lemma}\label{lem:Yamanouchi-permutation}
		Let $T$ be a standard Young tableau and let $\pi$ be its corresponding 
		Yamanouchi permutation.
		Then we have $\mathrm{Peak}(T) = \mathrm{Val}(\pi)$ and 
		$\mathrm{Val}(T) = \mathrm{Peak}(\pi)$.
	\end{lemma}

	We proceed to show that all the Yamanouchi permutations corresponding to 
	the SYT's of a given skew shape can be divided into several Knuth 
	equivalence classes.
	To this end, we need to define shape preserving transformations on 
	SYT's, which turn out to be equivalent to the Knuth transformations 
	on their corresponding Yamanouchi permutations.
	Hence we will also call the transformation the Knuth transformations 
	on SYT's.
	Let us first recall the Knuth transformations on permutations.
	
	\noindent {\bf The Knuth transformations on permutations} \\
	Let $\pi = \pi_1\pi_2\cdots \pi_n$ be a permutation in $\mathcal{S}_n$ 
	and let $a<b<c$ be three adjacent numbers $\pi_{i-1}\pi_i\pi_{i+1}$ in the 
	permutation $\pi$. 
	A {\em Knuth transformation} $\kappa_i$ of the permutation $\pi$ is its transformation
	into another permutation $\sigma$ that has one of the following four forms:
	
	\begin{enumerate}[label=(\roman*)]
		\item  $\pi =\cdots \; acb \;\cdots \; \longrightarrow \; \sigma=\cdots \; cab \;\cdots\;$;
		\item  $\pi =\cdots \; bca \;\cdots \; \longrightarrow \; \sigma=\cdots \; bac \;\cdots\;$;
		\item  $\pi =\cdots \; cab \;\cdots \; \longrightarrow \; \sigma=\cdots \; acb \;\cdots\;$;
		\item  $\pi =\cdots \; bac \;\cdots \; \longrightarrow \; \sigma=\cdots \; bca \;\cdots\;$.
		
	\end{enumerate}
	Thus each Knuth transformation switches two adjacent entries $a$ and $c$ provided an entry $b$ satisfying $a < b < c$ is located next to $a$ or $c$.
	From the definition of the transformations,  one can easily observe that 
	$\kappa_i$ is well-defined if and only if $i$ is 
	a peak or a valley of the permutation $\pi$ and $\kappa_i$ changes a peak to a valley, and vise versa.
	Two permutations $\pi$ and $\sigma$ are said to be {\em Knuth-equivalent},
	denoted by $\pi \stackrel{\text{K}}{\sim} \sigma$, if one of them can 
	be obtained from another by a sequence of Knuth transformations.
	Let $[\pi]$ denote the Knuth equivalence class that contains the permutation $\pi$.
	For example,  the five permutations as shown in Figure \ref{fig:Knuth} form a Knuth equivalence class, where the ones that differ by a single Knuth transformation are connected by an edge.
	
	
	\begin{figure}[H]
		\begin{center}
			\begin{tikzpicture}[font =\small , scale = 0.4, line width = 0.7pt]
				\tikzmath{\d = 6;};
				\foreach \x /\y in {1/32154,2/32514,3/35214,4/35241,5/32541}
				\node at (\d*\x,0) {$\y$};
				\foreach \i in {0,...,3} \draw (8+\i*\d,0)--(10+\i*\d,0);
			\end{tikzpicture}
		\end{center}
		\caption{An example of Knuth equivalence class of permutations}\label{fig:Knuth}
	\end{figure}
	
	\begin{theorem}(\cite{StanleyVol2},  Theorem A1.1.4 )\label{thm:Knuth-equi}
		Two permutations are Knuth-equivalent if and only if their insertion tableaux coincide.
	\end{theorem}
	
	Recall that a permutation $\pi$ has the same peak set and the same valley set with its recording tableau $Q$.
	Then the following lemma follows directly from the Theorem \ref{thm:Knuth-equi} and RSK algorithm.
	
	\begin{lemma}\label{lem:perm-SYT}
		Let $\pi$ be a permutation of length $n$.
		Then we have 
		$$\sum_{\tau \in [\pi]}t^{\mathrm{Peak}(\tau)} = \sum_{Q \in \mathrm{SYT}(\mathrm{sh}(\pi))}t^{\mathrm{Peak}(Q)}$$	and 
		$$\sum_{\tau \in [\pi]}t^{\mathrm{Val}(\tau)} = \sum_{Q \in \mathrm{SYT}(\mathrm{sh}(\pi))}t^{\mathrm{Val}(Q)}.$$
	\end{lemma}
	
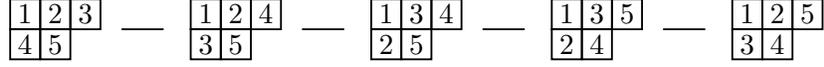
\begin{figure}[H]
		\begin{center}
			\begin{tikzpicture}[font =\small , scale = 0.4, line width = 0.7pt]
				\tikzmath{\d = 6;};
				\foreach \t in {0,...,4}
				{
					\foreach \i /\j in {3/1,2/2}
					{
						\foreach \k in {1,...,\i} \draw (\k+\t*\d,-\j)rectangle(\k +\t*\d+ 1,-\j -1);
					}
				}
				\foreach \i / \j / \k in {1/1/1,2/1/2,3/1/3,1/2/4,2/2/5} \node at (\i+0.5,-\j-0.5){$\k$};	
				
				\foreach \i / \j / \k in {1/1/1,2/1/2,3/1/4,1/2/3,2/2/5} \node at (\i+0.5+\d,-\j-0.5){$\k$};	
				
				\foreach \i / \j / \k in {1/1/1,2/1/3,3/1/4,1/2/2,2/2/5} \node at (\i+0.5+2*\d,-\j-0.5){$\k$};	
				
				\foreach \i / \j / \k in {1/1/1,2/1/3,3/1/5,1/2/2,2/2/4} \node at (\i+0.5+3*\d,-\j-0.5){$\k$};	
				
				\foreach \i / \j / \k in {1/1/1,2/1/2,3/1/5,1/2/3,2/2/4} \node at (\i+0.5+4*\d,-\j-0.5){$\k$};	
				
				\foreach \i in {0,...,3} \draw (4.7+\i*\d,-2)--(6.1+\i*\d,-2);
			\end{tikzpicture}
		\end{center}
		\caption{An example of Knuth equivalence class of SYT's}\label{fig:Knuth-SYT}
	\end{figure}

	\noindent {\bf The Knuth transformations on standard Young tableaux} \\
	Let $T$ be an SYT with $n \geq 3$ entries and let $i-1,i,i+1$ be three entries in $T$.
	A {\em Knuth transformation} $\kappa_i$ of the SYT $T$ is its transformation
	into another SYT $S$ that has one of the following four forms:
	
	\begin{enumerate}[label=(\roman*)]
		\item   If $i \in \mathrm{Val}(T)$ and $i+1$ appears in a lower row
		in $T$ than $i-1$, then $S$ is obtained from $T$ by switching $i-1$ and $i$ in $T$ and preserving all the other entries;
		\item If $i \in \mathrm{Val}(T)$ and $i+1$ does not appear in a lower row
		in $T$ than $i-1$, then $S$ is obtained from $T$ by switching $i$ and $i+1$ in $T$ and preserving all the other entries;
		\item  If $i \in \mathrm{Peak}(T)$ and $i+1$ does not appear in a lower row
		in $T$ than $i-1$, then $S$ is obtained from $T$ by switching $i-1$ and $i$ in $T$ and preserving all the  other entries;
		\item  If $i \in \mathrm{Peak}(T)$ and $i+1$ appears in a lower row
		in $T$ than $i-1$, then $S$ is obtained from $T$ by switching $i$ and $i+1$ in $T$ and preserving all the other entries.
		
	\end{enumerate}
	We remark that one can distinguish from the context whether $\kappa_i$ is defined on permutations or SYT's. 
	Two SYT's $T$ and $S$ are said to be {\em Knuth-equivalent},
	denoted by $T \stackrel{\text{K}}{\sim} S$, if one of them can 
	be obtained from another by a sequence of Knuth transformations.
	Let $[T]$ denote the Knuth equivalence class that contains the SYT $T$.
	For example,  the five SYT's as shown in Figure \ref{fig:Knuth-SYT} form a Knuth equivalence class, where the ones that differ by a single Knuth transformation are connected by an edge.
	
	The following basic facts are fairly straightforward
	and we omit the detailed  proofs here.
	
	\noindent {\bf Fact \#1}
	Let $T$ be an SYT.
	Then the two modified consecutive entries of $T$ under the Knuth transformation $\kappa_i$ are always not in the same row or in the same column, which implies that $\kappa_i(T)$ is an SYT of the same shape with $T$.
	
	\noindent {\bf Fact \#2}
	The Knuth transformation $\kappa_i$ is well-defined on $T$ if and only if $i$ is a peak or a valley of the SYT $T$ and $\kappa_i$ changes a peak to a valley, and vise versa.
	Moreover, $\kappa_i$ is an involution, namely, $\kappa_i^{2} = \mathrm{id}$, where $\mathrm{id}$
	is the identity transformation.
	
	\noindent {\bf Fact \#3}
	Let $T$ be an SYT of $n$ entries with $k$ and $k+1$ in different rows and columns and let 
	$S$ be the SYT obtained from $T$ by switching $k$ and 
	$k+1$ in $T$ and preserving other entries.
	Then we have that the Yamanouchi word $\alpha(S)$  of $S$ is the word obtained 
	from the Yamanouchi word $\alpha(T) = y_1y_2\cdots y_n$ of $T$ by switching 
	$y_k$ and $y_{k+1}$.
	Since $y_k \neq y_{k+1}$, we deduce that 
	the Yamanouchi permutation $\beta(S)$  of $S$ is the permutation obtained 
	from the Yamanouchi permutation $\beta(T) = \pi_1\pi_2\cdots \pi_n$ of $T$ by switching 
	$\pi_k$ and $\pi_{k+1}$.

	\begin{lemma}\label{lem:commutativity}
		Let $T$ be a standard Young tableau of skew shape $\lambda / \mu$ and 
		$i \in \mathrm{Peak}(T) \cup \mathrm{Val}(T)$.
		Then we have $$\beta\circ \kappa_i(T) = \kappa_i\circ \beta(T).$$
	\end{lemma}	  	
	
	\pf
	Let $\pi = \beta(T)$, $S = \kappa_i(T)$ and $\sigma = \kappa_i(\pi)$.
	For a better view, we
	outlines the sets and relationships as shown in Figure \ref{fig:struct1}.
	We need to show that $\sigma = \beta(S)$.
	We shall consider four cases according to the four forms of the 
	Knuth transformations on SYT's.
	
	\noindent {\bf Case (\rmnum{1}):}
	If $i \in \mathrm{Val}(T)$ and $i+1$ appears in a lower row
	in $T$ than $i-1$, then  $S$ is obtained from $T$ by switching $i-1$ and $i$ in $T$ and preserving all the other entries.
	Fact \#1 tells us that   $i-1$ and $i$ are in different rows and columns of $T$.
	Then by Fact \#3,  $\beta(S)$ is the permutation obtained from $\pi$ by switching $\pi_{i-1}$ and $\pi_i$ in $\pi$.
	On the other hand,  as $i \in \mathrm{Val}(T)$, then from Lemma \ref{lem:Yamanouchi-permutation},	we have $i \in \mathrm{Peak}(\pi)$.
	Note that $i+1$ appears in a lower row in $T$ than $i-1$.
	This implies that $\pi_{i-1} < \pi_{i+1}$.
	Thus $\pi_{i-1}\pi_i\pi_{i+1} = acb$ for some $a<b<c$.
	Then by the definition of $\kappa_i$, $\sigma$ is obtained from $\pi$ 
	by switching $\pi_{i-1}$ and $\pi_i$ in $\pi$,
	namely, $\sigma =\beta(S)$.
	
	\noindent {\bf Case (\rmnum{2}):}
	If $i \in \mathrm{Val}(T)$ and $i+1$ does not appear in a lower row
	in $T$ than $i-1$, then $S$ is obtained from $T$ by switching $i$ and $i+1$ in $T$ and preserving all the other entries.
	Fact \#1 tells us that   $i$ and $i+1$ are in different rows and columns of $T$.
	Then by Fact \#3,  $\beta(S)$ is the permutation obtained from $\pi$ by switching $\pi_{i}$ and $\pi_{i+1}$ in $\pi$.
	On the other hand,  as $i \in \mathrm{Val}(T)$, then from Lemma \ref{lem:Yamanouchi-permutation},	we have $i \in \mathrm{Peak}(\pi)$.
	Note that $i+1$ does not appear in a lower row in $T$ than $i-1$.
	This implies that $\pi_{i-1} > \pi_{i+1}$.
	Thus $\pi_{i-1}\pi_i\pi_{i+1} = bca$ for some $a<b<c$.
	Then by the definition of $\kappa_i$, $\sigma$ is obtained from $\pi$ 
	by switching $\pi_{i}$ and $\pi_{i+1}$ in $\pi$,
	namely, $\sigma = \beta(S)$.

	\noindent {\bf Case (\rmnum{3}):}
	If $i \in \mathrm{Peak}(T)$ and $i+1$ does not appear in a lower row
	in $T$ than $i-1$, then $S$ is obtained from $T$ by switching $i-1$ and $i$ in $T$ and preserving all the other entries.
	Fact \#1 tells us that   $i-1$ and $i$ are in different rows and columns of $T$.
	Then by Fact \#3,  $\beta(S)$ is the permutation obtained from $\pi$ by switching $\pi_{i-1}$ and $\pi_i$ in $\pi$.
	On the other hand,  as $i \in \mathrm{Peak}(T)$, then from Lemma \ref{lem:Yamanouchi-permutation},	we have $i \in \mathrm{Val}(\pi)$.
	Note that $i+1$ does not appear in a lower row in $T$ than $i-1$.
	We derive that $\pi_{i-1} > \pi_{i+1}$.
	Thus $\pi_{i-1}\pi_i\pi_{i+1} = cab$ for some $a<b<c$.
	Then by the definition of $\kappa_i$, $\sigma$ is obtained from $\pi$ 
	by switching $\pi_{i-1}$ and $\pi_i$ in $\pi$,
	namely, $\sigma = \beta(S)$.
	
	\noindent {\bf Case (\rmnum{4}):}
	If $i \in \mathrm{Peak}(T)$ and $i+1$ appears in a lower row
	in $T$ than $i-1$, then $S$ is obtained from $T$ by switching $i$ and $i+1$ in $T$ and preserving all the other entries.
	Fact \#1 tells us that   $i$ and $i+1$ are in different rows and columns of $T$.
	Then by Fact \#3,  $\beta(S)$ is the permutation obtained from $\pi$ by switching $\pi_{i}$ and $\pi_{i+1}$ in $\pi$.
	On the other hand,  as $i \in \mathrm{Peak}(T)$, then from Lemma \ref{lem:Yamanouchi-permutation},	we have $i \in \mathrm{Val}(\pi)$.
	Note that $i+1$ appears in a lower row in $T$ than $i-1$.
	This yields  that $\pi_{i-1} < \pi_{i+1}$.
	Thus $\pi_{i-1}\pi_i\pi_{i+1} = bac$ for some $a<b<c$.
	Then by the definition of $\kappa_i$, $\sigma$ is obtained from $\pi$ 
	by switching $\pi_{i}$ and $\pi_{i+1}$ in $\pi$,
	namely, $\sigma = \beta(S)$.

	So far, we have concluded that $\sigma = \beta(S)$ for any case, completing the proof.
	\qed

	\begin{figure}[H]
		\begin{center}
			\begin{tikzpicture}[font =\small , scale = 0.4, line width = 0.7pt]
				\tikzmath{\d = 4;}
				\node at (0,0) {$\pi$};
				\node at (0,\d) {$T$};
				\node at (\d,0) {$\sigma$};
				\node at (\d,\d) {$S$};
				\draw[<->](1,0)--(\d-1,0) node[below = 2mm,left]{$\kappa_i$};
				\draw[<-](0,1)--(0,\d-1) node[left = 2mm,below = 1mm]{$\beta$};
				\draw[<->](1,\d)--(\d-1,\d) node[above = 2mm,left]{$\kappa_i$};
				\draw[<-](\d,1)--(\d,\d-1) node[right = 2mm,below = 1mm]{$\beta$};			
			\end{tikzpicture}
		\end{center}
		\caption{The commutativity  of the maps $\beta$ and $\kappa_i$.} \label{fig:struct1}
	\end{figure}

	\begin{lemma}\label{lem:Knuth-equivalence}
		Let $T$ and $S$ be two standard Young tableaux of the shape $\lambda/\mu$ of size $n$ and 
		let $\pi$ and $\sigma$ be the corresponding Yamanouchi permutations of 
		$T$ and $S$, respectively.
		Then we have $$T \stackrel{\text{K}}{\sim} S \Leftrightarrow \pi \stackrel{\text{K}}{\sim} \sigma.$$
	\end{lemma}
	
	\pf
	If $T \stackrel{\text{K}}{\sim} S$, then there exists a sequence
	of Knuth transformations $\kappa_{i_1},\kappa_{i_2},\ldots,\kappa_{i_j}$
	such that $S = \kappa_{i_j}\circ\cdots \circ\kappa_{i_2}\circ\kappa_{i_1}(T)$.
	By Lemma \ref{lem:commutativity}, we have 
	$\beta(S) = \kappa_{i_j}\circ\cdots \circ\kappa_{i_2}\circ\kappa_{i_1}\circ\beta(T)$,
	that is, $\sigma = \kappa_{i_j}\circ\cdots \circ\kappa_{i_2}\circ\kappa_{i_1}(\pi)$.
	Hence $\pi \stackrel{\text{K}}{\sim} \sigma$.
	See Figure \ref{fig:struct2} for a better view.
	
	Conversely, if $\pi \stackrel{\text{K}}{\sim} \sigma$, 
	then there exists a sequence
	of Knuth transformations $\kappa_{i_1},\kappa_{i_2},\ldots,$ $\kappa_{i_j}$
	such that $\sigma = \kappa_{i_j}\circ\cdots  \circ\kappa_{i_2}\circ\kappa_{i_1}(\pi)$.
	Let $S' = \kappa_{i_j}\circ\cdots \circ\kappa_{i_2}\circ\kappa_{i_1}(T)$.
	Then we have $\beta(S') = \kappa_{i_j}\circ\cdots \circ\kappa_{i_2}\circ\kappa_{i_1}(\pi)=\sigma$.
	By the definition of the Knuth transformations on 
	SYT's, $S'$ is of shape $\lambda/\mu$.
	Since $\beta(S) = \sigma$ and $\beta$ is an injection from $\mathrm{SYT}(\lambda/\mu)$ to $\mathcal{S}_n$, 
	we have $S' = S$.
	Then $T \stackrel{\text{K}}{\sim} S$ follows directly.
	\qed
	
	\begin{figure}[H]
		\begin{center}
			\begin{tikzpicture}[font =\small , scale = 0.4, line width = 0.7pt]
				\tikzmath{\d = 4;}
				\node at (0,0) {$\pi$};
				\node at (\d,0) {$\pi^1$};
				\node at (2*\d,0) {$\pi^2$};
				\node at (2.45*\d,0) {$\cdots\cdots$};
				\node at (3*\d,0) {$\pi^{j-1}$};
				\node at (4*\d,0) {$\sigma$};
				
				\node at (0,\d) {$T$};
				\node at (\d,\d) {$T_1$};
				\node at (2*\d,\d) {$T_2$};
				\node at (2.45*\d,\d) {$\cdots\cdots$};
				\node at (3*\d,\d) {$T_{j-1}$};
				\node at (4*\d,\d) {$S$};
				
				\draw[<->](1,0)--(\d-1,0) node[below = 2mm,left]{$\kappa_{i_1}$};
				\draw[<->](1+\d,0)--(2*\d-1,0) node[below = 2mm,left]{$\kappa_{i_2}$};
				\draw[<->](1+3*\d,0)--(4*\d-1,0) node[below = 2mm,left]{$\kappa_{i_j}$};

				\draw[<-](0,1)--(0,\d-1) node[left = 2mm,below = 1mm]{$\beta$};
				\draw[<-](\d,1)--(\d,\d-1) node[left = 2mm,below = 1mm]{$\beta$};
				\draw[<-](2*\d,1)--(2*\d,\d-1) node[left = 2mm,below = 1mm]{$\beta$};
				\draw[<-](3*\d,1)--(3*\d,\d-1) node[left = 2mm,below = 1mm]{$\beta$};
				\draw[<-](4*\d,1)--(4*\d,\d-1) node[left = 2mm,below = 1mm]{$\beta$};
				
				\draw[<->](1,\d)--(\d-1,\d) node[above = 2mm,left]{$\kappa_{i_1}$};
				\draw[<->](1+\d,\d)--(2*\d-1,\d) node[above = 2mm,left]{$\kappa_{i_2}$};		
				\draw[<->](1+3*\d,\d)--(4*\d-1,\d) node[above = 2mm,left]{$\kappa_{i_j}$};
			\end{tikzpicture}
		\end{center}
		\caption{A better view for Lemma \ref{lem:Knuth-equivalence}.} \label{fig:struct2}
	\end{figure}

	Actually, the proof of Lemma \ref{lem:Knuth-equivalence} enables us to
	obtain the following lemma directly.
	
	\begin{lemma}\label{lem:equal}
		Let $T$ be a standard Young tableau of skew shape $\lambda / \mu$ and let $\pi$ be the Yamanouchi permutation of $T$.
		Then the map $\beta$ induces a bijection between the sets $[T]$ and $[\pi]$. 
	\end{lemma}

	The following lemma plays an 
	essential role in the proof of Theorem 
	\ref{thm:skewSYT}.
	
	\begin{lemma}\label{lem:1314}
		If $\mathrm{Peak}$ and $\mathrm{Val}$ are equidistributed over
		$\mathrm{SYT}(\nu)$ for any Young diagram $\nu$ of size $n$, 
		then  $\mathrm{Peak}$ and $\mathrm{Val}$ are also equidistributed over
		$\mathrm{SYT}(\lambda/\mu)$ for any skew diagram $\lambda/\mu$ of size $n$.
	\end{lemma}
	
	\pf
	Given any skew shape $\lambda/\mu$ of size $n$, 
	assume that $\mathrm{SYT}(\lambda /\mu)$ can be divided into $k$ 
	Knuth equivalence classes $[T_1],[T_2],\ldots,[T_k]$.
	Let $\pi^i$ be the corresponding Yamanouchi permutations of 
	$T_i$ and let $\upsilon^i$ be the shape of $\pi^i$ for 
	$1\leq i \leq k$.
	Then we have 
	\begin{align*}
		\sum_{T\in \mathrm{SYT}(\lambda/\mu)}t^{\mathrm{Peak}(T)} 
		&= \sum_{i=1}^k\sum_{T\in [T_i]}t^{\mathrm{Peak}(T)} \\
		&= \sum_{i=1}^k\sum_{\tau\in [\pi^i]}t^{\mathrm{Val}(\tau)} \quad   (\text{by Lemmas \ref{lem:Yamanouchi-permutation} and \ref{lem:equal} } )\\
		&= \sum_{i=1}^k\sum_{Q \in \mathrm{SYT}(\upsilon^i)}t^{\mathrm{Val}(Q)} \quad   (\text{by Lemma \ref{lem:perm-SYT}}   )\\
		&= \sum_{i=1}^k\sum_{Q \in \mathrm{SYT}(\upsilon^i)}t^{\mathrm{Peak}(Q)} \quad   (\text{by the hypothesis}  )\\
		&= \sum_{i=1}^k\sum_{\tau\in [\pi^i]}t^{\mathrm{Peak}(\tau)} \quad   (\text{by Lemma \ref{lem:perm-SYT}} )\\
		&= \sum_{i=1}^k\sum_{T\in [T_i]}t^{\mathrm{Val}(T)} \quad   (\text{by Lemmas \ref{lem:Yamanouchi-permutation} and \ref{lem:equal} } )\\
		&= \sum_{T\in \mathrm{SYT}(\lambda/\mu)}t^{\mathrm{Val}(T)}.
	\end{align*}
	This completes the proof.
	\qed

	For sets of positive integers, we define a {\em lexicographic order} on them, denoted by $\preceq$.
	Let $A = \{a_1,a_2,\ldots, a_k\}$ and $B = \{b_1,b_2,\ldots, b_r\}$ be two 
	sets of positive integers where $a_1<a_2<\cdots <a_k$ and $b_1<b_2<\cdots <b_r$.
	Throughout the paper, we always list a set of positive integers in increasing order.
	Define $A \preceq B$ if either $A = B$, or else for some $i$, 
	$$a_1 = b_1, a_2 = b_2,\ldots,a_i=b_i, a_{i+1} <b_{i+1},$$
	with the convention that $a_i = 0$ if $i > k$ and $b_i = 0$ if $i > r$.
	
	Given a skew diagram $\lambda/ \mu$,  define the multisets
	$$
	\mathcal{P}(\lambda/ \mu) = \{\mathrm{Peak}(T)\ \mid T \in  \mathrm{SYT}(\lambda/ \mu) \}
	$$
	and
	$$
	\mathcal{V}(\lambda/ \mu) = \{\mathrm{Val}(T)\ \mid T \in  \mathrm{SYT}(\lambda/ \mu)\}.
	$$
	For example,   
	we have $$\mathcal{P}((3,2,2)/(1,1)) = \mathcal{V}((3,2,2)/(1,1)) = \{\emptyset, \{2\}, \{2\}, \{2,4\}, \{2,4\}, \{3\}, \{3\}, \{3\}, \{3\}, \{4\}, \{4\}\},$$ 
	where the elements are listed in lexicographical order.
	See Table \ref{table:SYT} for details.

	\begin{table}[H]
		\centering
		\scriptsize
		\renewcommand\arraystretch{1.2}
		\caption{Equidistribution of $\mathrm{Peak}$ and $\mathrm{Val}$ over $\mathrm{SYT}((3,2,2)/(1,1))$.}\label{table:SYT}
		\vskip 2mm
		\begin{tabular}{|m{2cm}<{\centering}|c|c|m{2cm}<{\centering}|c|c|}
			\hline 
			$T$ &  $\mathrm{Peak}(T)$       &  $\mathrm{Val}(T)$     &    $T$   &   $\mathrm{Peak}(T)$     & $\mathrm{Val}(T)$    \\ \hline
			
			\vskip 1mm
			$\begin{array}[b]{*{3}c}\cline{2-3}
				&\lr{1}&\lr{3}\\\cline{2-3}
				&\lr{2}\\\cline{1-2}
				\lr{4}&\lr{5}\\\cline{1-2}
			\end{array}$ 
			&  $\{3\}$   & $\{2,4\}$     &  \vskip 1mm
			$\begin{array}[b]{*{3}c}\cline{2-3}
				&\lr{2}&\lr{3}\\\cline{2-3}
				&\lr{4}\\\cline{1-2}
				\lr{1}&\lr{5}\\\cline{1-2}
			\end{array}$  
			& $\{3\}$   & $\emptyset$                             \\ \hline
			
			\vskip 1mm
			$\begin{array}[b]{*{3}c}\cline{2-3}
				&\lr{1}&\lr{4}\\\cline{2-3}
				&\lr{2}\\\cline{1-2}
				\lr{3}&\lr{5}\\\cline{1-2}
			\end{array}$ 
			&  $\{4\}$   & $\{3\}$     &  \vskip 1mm
			$\begin{array}[b]{*{3}c}\cline{2-3}
				&\lr{2}&\lr{4}\\\cline{2-3}
				&\lr{3}\\\cline{1-2}
				\lr{1}&\lr{5}\\\cline{1-2}
			\end{array}$  
			& $\{2,4\}$   & $\{3\}$                             \\ \hline
			
			\vskip 1mm
			$\begin{array}[b]{*{3}c}\cline{2-3}
				&\lr{1}&\lr{5}\\\cline{2-3}
				&\lr{2}\\\cline{1-2}
				\lr{3}&\lr{4}\\\cline{1-2}
			\end{array}$ 
			&  $\emptyset$   & $\{3\}$     &  \vskip 1mm
			$\begin{array}[b]{*{3}c}\cline{2-3}
				&\lr{2}&\lr{5}\\\cline{2-3}
				&\lr{3}\\\cline{1-2}
				\lr{1}&\lr{4}\\\cline{1-2}
			\end{array}$  
			& $\{2\}$   & $\{4\}$                             \\ \hline

			\vskip 1mm
			$\begin{array}[b]{*{3}c}\cline{2-3}
				&\lr{1}&\lr{2}\\\cline{2-3}
				&\lr{4}\\\cline{1-2}
				\lr{3}&\lr{5}\\\cline{1-2}
			\end{array}$ 
			&  $\{2,4\}$   & $\{3\}$     &  \vskip 1mm
			$\begin{array}[b]{*{3}c}\cline{2-3}
				&\lr{1}&\lr{3}\\\cline{2-3}
				&\lr{4}\\\cline{1-2}
				\lr{2}&\lr{5}\\\cline{1-2}
			\end{array}$  
			& $\{3\}$   & $\{2\}$                             \\ \hline
			
			\vskip 1mm
			$\begin{array}[b]{*{3}c}\cline{2-3}
				&\lr{1}&\lr{2}\\\cline{2-3}
				&\lr{3}\\\cline{1-2}
				\lr{4}&\lr{5}\\\cline{1-2}
			\end{array}$ 
			&  $\{2\}$   & $\{4\}$     &  \vskip 1mm
			$\begin{array}[b]{*{3}c}\cline{2-3}
				&\lr{1}&\lr{4}\\\cline{2-3}
				&\lr{3}\\\cline{1-2}
				\lr{2}&\lr{5}\\\cline{1-2}
			\end{array}$  
			& $\{4\}$   & $\{2\}$  	                              \\ \hline
			
			\vskip 1mm
			$\begin{array}[b]{*{3}c}\cline{2-3}
				&\lr{1}&\lr{5}\\\cline{2-3}
				&\lr{3}\\\cline{1-2}
				\lr{2}&\lr{4}\\\cline{1-2}
			\end{array}$
			& $\{3\}$   & $\{2,4\}$ 
			&  	   	   &      &                     \\ \hline
		\end{tabular}
	\end{table}
	
	Let $M = \{M_1,M_2,\ldots,M_k\}$ be a multiset in which each element $M_i$ is a set of positive integers and let $A$ be a set of positive integers.
	Define $$M[A] = \{M_1 \cap A, M_2\cap A, \ldots, M_k \cap A \}.$$
	We call $M[A]$ the restriction of $M$ into $A$.
	For example,  we have 
	$$\mathcal{P}((3,2,2)/(1,1)) [\{2,3\}] = \{\emptyset, \emptyset, \emptyset,  \{2\},\{2\},  \{2\},  \{2\}, \{3\},\{3\},  \{3\}, \{3\} \}.$$
	For positive integers $a,b$ with $a\leq b$, we denote by $[a,b]$ the set of all integers $j$ with $a\leq j\leq b$. 
	For a set $A$ of positive integers and an positive integer $k$, we denote by $A+k$ the set obtained 
	from $A$ by adding $k$ to each element of $A$.
	Similarly, for an SYT $T$ and a positive integer $k$, we denote by $T+k$ the tableau obtained from $T$ by 
	increasing each entry in $T$ by $k$.	
	The following lemma will play  an essential role in the proof of Theorem  \ref{thm:skewSYT}. 
	
	\begin{lemma}\label{lem:splice}
		Let $n$ be a positive integer and let $A$ be any subset of $[n]$ such that  $k,k+1 \notin A$ for  some positive integer $k$ $(1\leq k \leq n-1)$.		
		If   $\mathrm{Peak}$ and $\mathrm{Val}$ are  equidistributed over $\mathrm{SYT}(\lambda' / \mu')$ for  any skew diagram $\lambda'/ \mu'$ of size less than $n$, 
		then for any skew  diagram $\lambda/ \mu$ of size $n$, we have 
		$$\mathcal{P}( \lambda/ \mu )[A] = \mathcal{V}( \lambda/ \mu)[A].$$
	\end{lemma}
	
	\pf
	Let $B = [n] \setminus [k, k+1]$. 
	Note that $k,k+1 \notin A$, we have $A \subseteq B$.
	One can easily verify that $M[A] = M[B][A]$ for any multiset $M$.
	Hence it is sufficient to prove $\mathcal{P}( \lambda/ \mu )[B] = \mathcal{V}( \lambda/\mu )[B]$.
	
\begin{figure}[H]
		\begin{center}
			\begin{tikzpicture}[font =\small , scale = 0.3, line width = 0.7pt]
				\tikzmath{\r = 13;}
				\node at (4,5) {\footnotesize{$P$}};
				\node at (3,0.5) {\footnotesize{$Q+k$}};
				\node at (4,-3) {$T$};			
				\draw[->] (9,3)--(12,3) node[left = 10,above]{$\omega$};
				\foreach \d in {0,\r}
				{
					\draw (0+\d,-1)--(6+\d,-1)--(6+\d,2)--(8+\d,2)--(8+\d,5)--(10+\d,5)--(10+\d,8)--(4+\d,8);
					\draw (0+\d,-1)--(0+\d,4)--(2+\d,4)--(2+\d,6)--(4+\d,6)--(4+\d,8);
					\draw 		(0+\d,2)--(4+\d,2)--(4+\d,4)--(6+\d,4)--(6+\d,6)--(8+\d,6)--(8+\d,8);
				}
				\node at (4+\r,5) {\footnotesize{$\omega_1(P)$}};
				\node at (3+\r,0.5) {\footnotesize{$\omega_2(Q)+k$}};
				\node at (4+\r,-3) {$\omega(T)$};
			\end{tikzpicture}
		\end{center}
		\caption{A better view for the bijection $\omega$.} \label{fig:splice}
	\end{figure}
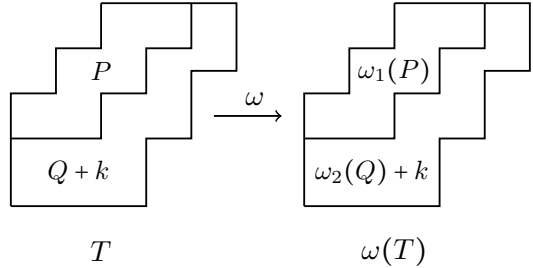

	We proceed to construct a bijection $\omega: \mathrm{SYT}(\lambda/ \mu) \rightarrow \mathrm{SYT}(\lambda/ \mu)$.
	For an SYT $T\in \mathrm{SYT}(\lambda/ \mu)$, 
	let $P$ be the SYT which is obtained from $T$ by reading the 
	entries in $[k]$. 
	And let $Q$ be the SYT which is obtained from $T$ by reading the 
	entries in $[k+1,n]$ and decreasing each entry by $k$.
	Assume that 
	$Q$ is of the skew shape $\lambda /\mu'$. Then the shape of $P$ is given by $\mu'/ \mu$.  
	By the hypothesis,  there exist a shape preserving 
	bijection, say $\omega_1: \mathrm{SYT}(\mu'/ \mu) \rightarrow \mathrm{SYT}(\mu'/ \mu)$,
	such that  $\mathrm{Peak}(R) = \mathrm{Val}(\omega_1(R))$ for any $R\in \mathrm{SYT}(\mu'/ \mu)$ 
	and a shape preserving  bijection, say $\omega_2: \mathrm{SYT}(\lambda /\mu') \rightarrow \mathrm{SYT}(\lambda /\mu')$, such that $\mathrm{Peak}(R) = \mathrm{Val}(\omega_2(R))$ for any $R\in \mathrm{SYT}(\lambda /\mu')$.
	Then $\omega(T)$ is defined to be the SYT by merging $\omega_1(P)$ and $ \omega_2(Q)+k$.
	See Figure \ref{fig:splice} for an illustration.
	Since both $\omega_1$ and $\omega_2$ are shape preserving bijections, 
	$\omega$ is a bijection.
	It is easily checked that $\mathrm{Peak}(T) \cap B = \mathrm{Peak}(P)\cup (\mathrm{Peak}(Q) +k)  = \mathrm{Val}(\omega_1(P))\cup (\mathrm{Val}(\omega_2(Q)) +k ) = \mathrm{Val}(\omega(T))\cap B$.
	To conclude, we construct a bijection $\omega$ from $\mathrm{SYT}(\lambda/ \mu)$ to 
	itself such that $\mathrm{Peak}(T) \cap B = \mathrm{Val}(\omega(T))\cap B$ for 
	any $T\in \mathrm{SYT}(\lambda/\mu)$.
	Thus $\mathcal{P}( \lambda/\mu )[B] = \mathcal{V}( \lambda/ \mu)[B]$, completing the proof.
	\qed

	Let $\lambda = (\lambda_1,\lambda_2,\ldots, \lambda_k)$ be an integer partition.
	The {\em rank} of $\lambda$, denoted by $\mathrm{rank}(\lambda)$, is 
	defined to be the largest  $i$ for which $\lambda_i \geq i$. 
	Equivalently, $\mathrm{rank}(\lambda)$ is the length of the main diagonal in the Young diagram of $\lambda$. 
	
	\begin{lemma}\label{lem:empty}
		Let $\lambda$ be a Young diagram of size $n$.
		\begin{enumerate}[label=\upshape(\roman*)]
		\item	If $\mathrm{rank}(\lambda)\geq 2$, then there  does not exist  $T \in \mathrm{SYT}(\lambda)$ with $\mathrm{Peak}(T) = \emptyset$ and there does not exist $T \in \mathrm{SYT}(\lambda)$ with $\mathrm{Val}(T) = \emptyset$.
		\item 			If $\mathrm{rank}(\lambda) = 1$, then there exists exactly one $T \in \mathrm{SYT}(\lambda)$ with $\mathrm{Peak}(T) = \emptyset$ and there exists exactly one  $T \in \mathrm{SYT}(\lambda)$ with $\mathrm{Val}(T)
			 = \emptyset$.
		\end{enumerate}
	\end{lemma}
	
	\pf
	Assume that $\mathrm{rank}(\lambda)\geq 2$. 
	We proceed to show that there does not exist $T \in \mathrm{SYT}(\lambda)$ with $\mathrm{Peak}(T) = \emptyset$.
	If not, let $T$ be the SYT in $\mathrm{SYT}(\lambda)$ such that
	$\mathrm{Peak}(T) = \emptyset$ and let $y = y_1y_2\cdots y_n$ be the corresponding Yamanouchi word of $T$.
	Then by Lemma \ref{lem:Yamanouchi-word}, we have $y \in \mathcal{Y}(\lambda)$ with $\mathrm{Val}(y) = \emptyset$.
	Since  $\mathrm{rank}(\lambda)\geq 2$, there exist at least two occurrences of 1 and 
	at least two occurrences of 2 in $y$.
	If follows that there exists some $i$ such that $y_i \geq y_{i+1}$.
	Let $i$ be the smallest such integer.
	As $y \in \mathcal{Y}(\lambda)$, we deduce that $y_j = j$ for  all $1\leq j\leq i$ and $y_{i+1} = 1$.
	Note that $y$ contains at least two occurrences of 2.
	It yields that there exists some $k\geq i+1$ such that $y_k<y_{k+1}$.
	Let $k$ be the smallest such integer.
	Then we have $y_{k-1}\geq y_k <y_{k+1}$, namely, $k \in \mathrm{Val}(y)$,
	a contradiction with the fact $\mathrm{Val}(y) = \emptyset$.
	Hence  there does not exist $T \in \mathrm{SYT}(\lambda)$ with $\mathrm{Peak}(T) = \emptyset$.
	
	Now we prove that there does not exist $T \in \mathrm{SYT}(\lambda)$ with $\mathrm{Val}(T) = \emptyset$.
	If not, let $T$ be the SYT in $\mathrm{SYT}(\lambda)$ such that
	$\mathrm{Val}(T) = \emptyset$ and let $y = y_1y_2\cdots y_n$ be the corresponding Yamanouchi word of $T$.
	Then by Lemma \ref{lem:Yamanouchi-word}, we have $y \in \mathcal{Y}(\lambda)$ with $\mathrm{Peak}(y) = \emptyset$.
	Note that $y_1= 1$ and $y$ contains at least two occurrences of 2.
	If follows that there exists some $i$ such that $y_i < y_{i+1}$.
	Let $i$ be the smallest such integer.
	As $y \in \mathcal{Y}(\lambda)$, we deduce that $y_j= 1$ for all $1\leq j\leq i$	and $y_{i+1} = 2$.
	Since $y$ contains at least two occurrences of 2,  there exists some $k\geq i+1$ such that $y_k \geq y_{k+1}$.
	Let $k$ be the smallest such integer.
	Then we have  $y_{k-1}< y_k \geq y_{k+1}$, namely, $k \in \mathrm{Peak}(y)$,
	a contradiction with the fact $\mathrm{Peak}(y) = \emptyset$.
	Hence  there does not exist $T \in \mathrm{SYT}(\lambda)$ with $\mathrm{Val}(T) = \emptyset$, completing the proof of 
	(\upshape \rmnum{1}).
	
	Now we proceed to prove (\upshape \rmnum{2}).
	Assume that $\mathrm{rank}(\lambda) = 1$.
	By similar arguments as the proof of (\upshape \rmnum{1}),
	one can easily derive that there exists exactly one Yamanouchi word of the form $y = 12\cdots k 11\cdots 1   \in \mathcal{Y}(\lambda)$ with $\mathrm{Val}(y) = \emptyset$ and there exists exactly one Yamanouchi word of the form $y = 11\cdots 12\cdots k   \in \mathcal{Y}(\lambda)$ with $\mathrm{Peak}(y) = \emptyset$.
	Then (\upshape \rmnum{2}) follows directly from Lemma \ref{lem:Yamanouchi-word}.
	\qed

	Now we are ready for the proof of Theorem \ref{thm:skewSYT}.
	
	\noindent{\bf  Proof of Theorem \ref{thm:skewSYT}.}
	By Lemma \ref{lem:1314}, it is sufficient to prove that
	 $\mathrm{Peak}$  and $\mathrm{Val}$  are equidistributed over
	$\mathrm{SYT}(\lambda)$ for any Young diagram $\lambda$ of size $n$.
	We prove the assertion   by induction on $n$.
	It is routine to check the assertion  for $n\leq 3$.
	Assume that the assertion holds for all $m < n$ with $n \geq 4$.
	We proceed to prove $\mathcal{P}(  \lambda) = \mathcal{V}( \lambda)$
	for any Young diagram $\lambda$ of size $n$.

	If we let $A = [1,n-2]$, then by the induction hypothesis,  Lemma \ref{lem:1314} and Lemma \ref{lem:splice},
	we have $\mathcal{P}(  \lambda)[A] = \mathcal{V}(  \lambda)[A]$.
	Let $\mathcal{P}=\mathcal{P}(\lambda) = \{P_1,P_2,\ldots, P_{t}\}$ where
	$t = f^{\lambda}$ and 
	$P_i$ ($1\leq i \leq t$) are sorted according to the following rules:
	\begin{itemize}
		\item $P_1\cap A \preceq P_2\cap A \preceq\cdots \preceq P_{t}\cap A$;
		\item for $j>i$, if $P_i \cap A = P_{j}\cap A$ and $n-1 \in P_i$, then $n-1 \in P_{j}$.
	\end{itemize} 
	Similarly, let $\mathcal{V}=\mathcal{V}(\lambda) = \{V_1,V_2,\ldots, V_{t}\}$ where
	$V_i$ ($1\leq i \leq t$) are sorted according to the following rules:
	\begin{itemize}
		\item $V_1\cap A \preceq V_2\cap A \preceq\cdots \preceq V_{t}\cap A$;
		\item for $j>i$, if $V_i \cap A = V_{j}\cap A$ and $n-1 \in V_i$, then $n-1 \in V_{j}$.
	\end{itemize} 
	It follows that $P_i\cap A = V_i\cap A$ for $1\leq i\leq t$.
	Combining the fact that $P_i, V_i \subseteq [2,n-1]$, 
	we have either $P_i = V_i$ or $P_i$ and $V_i$ differ by one element $n-1$.
	
	We shall prove that $P_i = V_i$ for $1\leq i\leq t$.
	Assume on the contrary  there exists some $i$ such that $P_i \neq V_i$.
	Choose $i$ to be the smallest such integer.
	Since $P_i$ and $V_i$   differ by one element $n-1$, without loss of generality, 
	assume that $P_i = V_i\cup  \{n-1\}$ and $n-1 \notin V_i$.
	
	\noindent{\bf Claim 1.} $V_i \neq  \emptyset$.\\
	If not,  then we have  $V_i = \emptyset$ and $P_i = \{n-1\}$.
	Notice that $P_i\cap A = \emptyset$.
	By the ordering rules of the elements in $\mathcal{P}$, it is easily checked that 
	$P_j \neq \emptyset$ for all $j \geq i$.
	Note that 
	$P_j=V_j$ for all $j<i$. 
	One can easily check that $\mathcal{V}$ contains more empty sets than $\mathcal{P}$,
	which contradicts Lemma \ref{lem:empty}.
	Hence we have $V_i \neq \emptyset$.

	
	By Claim 1, we assume that $V_i=\{b_1, b_2, \ldots, b_k\}$ with $2\leq b_1<b_2<\ldots <b_k< n-1$. 
	\noindent{\bf Claim 2.} $b_1 = 2$.\\
	If not,  then	$P_i\cap B = \{b_1, b_2, \ldots, b_k,n-1\}$ and 
	$V_i\cap B = \{b_1, b_2, \ldots, b_k\}$, where $B=[3, n-1]$.
	Again by the induction hypothesis,  Lemma \ref{lem:1314} and Lemma \ref{lem:splice}, we have $\mathcal{P}[B] = \mathcal{V}[B]$.
	We assert that $P_j\cap B \neq \{b_1, b_2, \ldots, b_k\}$ for $j > i$.
	If not, we have either $P_j =   \{b_1, b_2, \ldots, b_k\}$ or 
	$P_j =   \{2,b_1, b_2, \ldots, b_k\}$.
	In both cases, we have $P_j \cap A = P_j$.
	It yields a contradiction with 
	the ordering rules of the elements in $\mathcal{P}$ as $P_i \cap A = \{b_1, b_2, \ldots, b_k\}$ and $b_1 \geq 3$.
	Hence the assertion holds.
	Recall that 
	$P_j=V_j$ for all $j<i$. 
	Hence  $\mathcal{V}[B]$ contains more $\{b_1, b_2, \ldots, b_k\}$'s than $\mathcal{P}[B]$,
	a contradiction. 
	 Hence, the claim is proved.

	By Claim 2, there exists at least one integer $\ell$ ($1\leq \ell\leq k$) such that $b_j=2j$ for all $j\leq \ell$. Choose $\ell$ to be the largest such integer.

	\noindent{\bf Claim 3.} $\ell<k$.\\
	If not, we have $\ell = k$ and $b_k = 2k$, implying that 
	$P_i  = \{2,4,\dots,2k,n-1\}$ and 
	$V_i  = \{2,4,\dots,2k\}$.
	Let $C = [2,n-1]\setminus[2k-1,2k]$. 
	Again by the induction hypothesis,  Lemma \ref{lem:1314} and Lemma \ref{lem:splice},
	we have 
	$\mathcal{P}[C] = \mathcal{V}[C]$. 
	We consider two cases for $k$.\\
	\noindent{\bf Case 1.} $k>1$.\\
	It is easily seen that   
	$P_i \cap C = \{2,4,\dots,2k-2,n-1\}$   and 
	$V_i \cap C = \{2,4,\dots,2k-2\}$.
   We assert that $P_j\cap C \neq \{2,4,\dots,2k-2\}$ for $j > i$.
    If not, we have either $P_j =   \{2,4,\dots,2k-2\}$ or 
    $P_j =   \{2,4,\dots,2k-2,2k\}$.
    For both cases of $P_j$,  it will lead to a contradiction with 
    the ordering rules of the elements in $\mathcal{P}$.
    Hence the assertion holds.
     Recall that 
	$P_j=V_j$ for all $j<i$. 
	Hence  $\mathcal{V}[C]$ contains more $\{2,4,\dots,2k-2\}$'s than $\mathcal{P}[C]$,
	a contradiction.  \\
	\noindent{\bf Case 2.} $k=1$.\\
	Then we have $P_i  =\{2,n-1\}$ and $V_i = \{2\}$.
	 It yields that $P_i \cap C = \{n-1\}$  and $V_i \cap C =  \emptyset$.
	 Similarly, by the ordering rules of of the elements in $\mathcal{P}$,
	 we derive that $P_j\cap C \neq \emptyset$ for $j > i$. Recall that 
	$P_j=V_j$ for all $j<i$. 
	Hence  $\mathcal{V}[C]$ contains more empty sets than $\mathcal{P}[C]$,
	a contradiction.
	Hence the claim is proved.

	\noindent{\bf Claim 4.} $b_{\ell+1} = 2\ell+2$.\\
	If not, we have $b_{\ell+1} > 2\ell+2$ and $n -1 > 2\ell +4$.
	Let  $D = [2,n-1] \setminus [2\ell+1,2\ell+2]$. 
		Again by the induction hypothesis,  Lemma \ref{lem:1314} and Lemma \ref{lem:splice}, we have 
	$\mathcal{P}[D] = \mathcal{V}[D]$.
	 One can easily check that
	$P_i \cap D = P_i = \{2,4,\dots,2\ell,b_{l+1},\ldots, b_k,n-1\}$ and 
	$V_i \cap D = V_i = \{2,4,\dots,2\ell,b_{l+1},\ldots, b_k\}$.
Again by the ordering rules of the elements in $\mathcal{P}$, 
 we derive that $P_j\cap D \neq  V_i = \{2,4,\dots,2\ell,b_{l+1},\ldots, b_k\}$ for $j > i$.  Recall that 
	$P_j=V_j$ for all $j<i$. 
	This implies that $\mathcal{V}[D]$ contains more $ \{2,4,\ldots,2\ell,b_{l+1},\cdots, b_k\}$'s than $\mathcal{P}[D]$, a contradiction. Hence, we have $b_{\ell+1} = 2\ell+2$.
	
	Combining Claims 3 and 4, we have $b_{j} = 2j$ for all $j\leq \ell+1$ and $\ell+1\leq k$. This yields a contradiction with the choice of $\ell$.  Hence, we have concluded that  $P_i = V_i$ for $1\leq i\leq t$ as desired, completing the proof. 
	\qed

	We remark that the joint distribution of $\mathrm{Peak}$ and $\mathrm{Val}$ is in general not symmetric over $\mathrm{SYT}(\lambda /\mu)$.  For example,   we have 
	$$\sum_{T \in \mathrm{SYT}((3,2,2)/(1,1)) }t^{\mathrm{Peak}(T)}q^{\mathrm{Val}(T)} = 2t_3q_2q_4+2t_2t_4q_3+t_4q_3+q_3+t_3+2t_2q_4+t_3q_2+t_4q_2.
	$$

	The following theorem follows directly from Theorem \ref{thm:skewSYT} and Lemma \ref{lem:Yamanouchi-word}. 
	
	\begin{theorem}\label{thm:y}
		For a   skew diagram $\lambda /\mu$  of size $n$,  
		 $\mathrm{Peak}$  and $\mathrm{Val}$  are equidistributed over
		$\mathcal{Y}(\lambda /\mu)$,  that is,   
		$$\sum_{y \in \mathcal{Y}(\lambda /\mu) }t^{\mathrm{Peak}(y)} = \sum_{y \in \mathcal{Y}(\lambda /\mu)  }t^{\mathrm{Val}(y)}.$$
	\end{theorem}
	
	

	We conclude this section with  the  equidistribution of the peak set and the valley set on permutations and involutions with given shape. 
	The following theorem   follows directly from 
 	Theorem \ref{thm:skewSYT} and	Lemma \ref{lem:perm-SYT}.
	
	\begin{theorem}\label{thm:Knuth-peak-valley}
		Let $\pi$ be a permutation of length $n$. 
		Then $\mathrm{Peak}$ and $\mathrm{Val}$ are equidistributed over
		$[\pi]$, namely, 
		$$\sum_{\tau \in [\pi] }t^{\mathrm{Peak}(\tau)} = \sum_{\tau  \in [\pi] }t^{\mathrm{Val}(\tau)}.$$
	\end{theorem}
	
	From Theorem \ref{thm:Knuth-equi}, 
	we deduce that the set of permutations with a given shape $\lambda$ can be 
	divided into $f^{\lambda}$ Knuth equivalence classes,
	in which any two different Knuth equivalence classes have the same distribution of
	peak set (or valley set).
	Then the following theorem follows directly.
	
	\begin{theorem}\label{thm:Shape-peak-valley}
		Let $\lambda$ be a Young diagram of size $n$ and let $\pi$ be any permutation of shape $\lambda$.
		Then $\mathrm{Peak}$ and $\mathrm{Val}$ are equidistributed over the set of permutations
		of  shape $\lambda$, namely, 
		$$\sum_{\tau\in \mathcal{S}_n: \mathrm{sh}(\tau) = \lambda }t^{\mathrm{Peak}(\tau)} = \sum_{\tau\in \mathcal{S}_n: \mathrm{sh}(\tau) = \lambda  }t^{\mathrm{Val}(\tau)} = f^{\lambda}  \sum_{\tau \in [\pi] }t^{\mathrm{Peak}(\tau)} = f^{\lambda}\sum_{\tau  \in [\pi] }t^{\mathrm{Val}(\tau)}.$$
	\end{theorem}
	
	Notice that  $\pi \xrightarrow{\mathrm{RSK}} (P,Q)$ if only if $\pi^{-1}\xrightarrow{\mathrm{RSK}} (Q,P)$.
	It follows that RSK algorithm associates an involution $\pi$
	with an SYT $T$ such that $\mathrm{Peak}(\pi) = \mathrm{Peak}(T)$ and $\mathrm{Val}(\pi) = \mathrm{Val}(T)$.
	Hence the following theorem follows directly from Theorem \ref{thm:skewSYT}.
	
	\begin{theorem}\label{thm:Involution-peak-valley}
		Let $\lambda$ be a Young diagram of size $n$.
		Then $\mathrm{Peak}$ and $\mathrm{Val}$ are equidistributed over the set of involutions  of shape $\lambda$, namely, 
		$$\sum_{\tau\in \mathcal{I}_n: \mathrm{sh}(\tau) = \lambda }t^{\mathrm{Peak}(\tau)} = \sum_{\tau\in \mathcal{I}_n: \mathrm{sh}(\tau) = \lambda  }t^{\mathrm{Val}(\tau)}.$$
	\end{theorem}

	\section{Peaks and valleys on transversals}\label{sec:3}
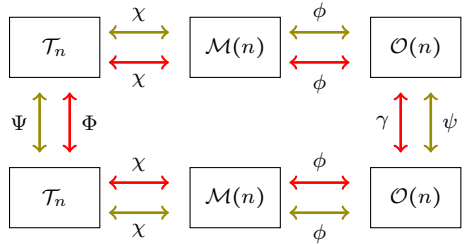
\begin{figure}[H]
		\centering
		\begin{tikzpicture}[font = \scriptsize , scale = 0.4]
			\tikzmath{\w = 6;}
			\tikzmath{\h = 5;}
			\foreach \x/\y/\z in {0/0/$\mathcal{T}_n$,\w/0/$\mathcal{M}(n)$,2*\w/0/$\mathcal{O}(n)$,0/\h/$\mathcal{T}_n$,\w/\h/$\mathcal{M}(n)$,2*\w/\h/$\mathcal{O}(n)$}
			{
				\draw (\x,\y)rectangle(\x+3,\y+2);
				\node at (\x+1.5,\y+1){\z};
			}
			\draw[<->,color = olive,line width=1](1,2.5)--(1,4.5) node[color = black,above=-4mm,left]{$\Psi$};
			\draw[<->,color = red,line width=1](1+2*\w,2.5)--(1+2*\w,4.5) node[color = black,above=-4mm,left]{$\gamma$};
			
			\draw[<->,color = red,line width=1](2,2.5)--(2,4.5) node[color = black,above=-4mm,right]{$\Phi$};
			\draw[<->,color = olive,line width=1](2+2*\w,2.5)--(2+2*\w,4.5) node[color = black,above=-4mm,right]{$\psi$};
			
			\draw[<->,color = olive,line width=1](3.3,0.5)--(5.3,0.5) node[color = black,left=4mm,below]{$\chi$};
			\draw[<->,color = red,line width=1](3.3,1.5)--(5.3,1.5) node[color = black,left=4mm,above]{$\chi$};
			
			\draw[<->,color = red,line width=1](3.3,0.5+\h)--(5.3,0.5+\h) node[color = black,left=4mm,below]{$\chi$};
			\draw[<->,color = olive,line width=1](3.3,1.5+\h)--(5.3,1.5+\h) node[color = black,left=4mm,above]{$\chi$};
			
			\draw[<->,color = olive,line width=1](3.3+\w,0.5)--(5.3+\w,0.5) node[color = black,left=4mm,below]{$\phi$};
			\draw[<->,color = red,line width=1](3.3+\w,1.5)--(5.3+\w,1.5) node[color = black,left=4mm,above]{$\phi$};
			
			\draw[<->,color = red,line width=1](3.3+\w,0.5+\h)--(5.3+\w,0.5+\h) node[color = black,left=4mm,below]{$\phi$};
			\draw[<->,color = olive,line width=1](3.3+\w,1.5+\h)--(5.3+\w,1.5+\h) node[color = black,left=4mm,above]{$\phi$};
			
		\end{tikzpicture}
		\caption{A diagrammatic summary of the sets and bijections.} \label{fig:struct}
	\end{figure}

	In this section, we aim to investigate the distribution of the  set-valued statistics Peak and Val on   transversals, thereby proving Theorems \ref{thm:transversal}, \ref{thm:transversal-general},
	\ref{thm:symm-transversal}, \ref{thm:symm-transversal-general} and \ref{thm:AI}.  To this end, 
	we shall establish an involution $\Phi$ on  $\mathcal{T}_n$ and a bijection $\Psi$
	from $\mathcal{T}_n$ to itself via matchings and oscillating tableaux  as intermediate structures. 
	The involution $\Phi$ is  achieved by employing    the conjugate map $\gamma$ on  $\mathcal{O}(n)$, 
	while the bijection $\Psi$ is  accomplished by constructing  a bijection $\psi$ from  $\mathcal{O}(n)$
	to itself based on Theorem \ref{thm:skewSYT}.  Figure \ref{fig:struct} outlines our sets and maps of interest.

	\subsection{The bijection $\chi$ from transversals to  matchings}
	
\begin{figure}[H]
		\begin{center}
			\begin{tikzpicture}[font =\small , scale = 0.5, line width = 0.7pt]
				\foreach \x in {1,...,6}
				{
					\filldraw[black](\x,0)circle(2pt);
					\node at(\x,-0.5) {\x};
				}
				\draw (1,0)..controls(1.66,1.5)and(2.32,1.5)..(3,0);
				\draw (2,0)..controls(3,2)and(4,2)..(5,0);
				\draw (4,0)..controls(4.66,1.5)and(5.32,1.5)..(6,0);				
			\end{tikzpicture}
		\end{center}
		\caption{The linear representation of the matching $\{(1,3),(2,5),(4,6)\}$.}\label{fig:matching}
	\end{figure}
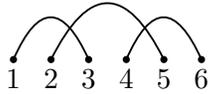

	A (perfect) {\em matching} $M$ of $[2n]$ is a collection of $n$ pairs 
	$\{(i_1,j_1),(i_2,j_2),\ldots,(i_n,j_n)\}$ with $i_k < j_k$ ($1\leq k \leq n$) such that each number of $[2n]$ appears exactly once.
	A pair $(i_k,j_k)$ ($1\leq k \leq n$) is called an {\em arc} of $M$,
	where $i_k$ is called an {\em opener} while $j_k$ is called a {\em closer}.
	Let $\mathcal{M}(n)$ denote the set of matchings of $[2n]$.
	Given a matching $M$ of $[2n]$, it can be represented by a graph $G$ with the vertex set $[2n]$ whose edge set consists of arcs of $M$.
	We usually draw the vertices of $G$ on a horizontal line in increasing order and draw the arcs of $G$ above the horizontal line.
	Such a graph is called the {\em linear representation} of the matching $M$.
	For example, Figure \ref{fig:matching} is the linear representation of the matching $\{(1,3),(2,5),(4,6)\}$.

	A {\em $k$-crossing} of a matching $M$ is a $k$-subset
	$(i_1, j_1), (i_2, j_2), \ldots, (i_k, j_k)$ of the arcs of $M$
	such that $i_1<i_2<\cdots <i_k<j_1<j_2<\cdots<j_k$.
	Similarly, a {\em $k$-nesting} is a $k$-subset
	$(i_1, j_1), (i_2, j_2), \ldots, (i_k, j_k)$ of the arcs of $M$
	such that $i_1<i_2<\cdots <i_k<j_k<\cdots<j_2<j_1$.
	A matching without any $k$-crossing is called {\em $k$-noncrossing matching}
	and a matching without any $k$-nesting is called {\em $k$-nonnesting matching}.
	Let $\mathcal{CM}_k(n)$ and $\mathcal{NM}_k(n)$ denote the sets of $k$-noncrossing and $k$-nonnesting matchings of $[2n]$, respectively.
	
	Given a matching $M = \{(i_1,j_1),(i_2,j_2),\ldots,(i_n,j_n)\}$ of $[2n]$, 
	define $M^r = \{(2n+1- j_1,2n+1-i_1),(2n+1 -j_2,2n+1-i_2),\ldots,(2n+1-j_n,2n+1-i_n)\}$.
	We call $M^r$ the {\em reverse} of the matching $M$.
	If $M = M^r$, we say $M$ is {\em bilaterally symmetric}.
	It can be easily seen that $M$ is bilaterally symmetric
	if and only if its linear representation is symmetric along the vertical line $x={2n+1\over 2}$.
	Figure \ref{fig:matching} illustrates a bilaterally symmetric matching.
	Let $\mathcal{SM}(n)$ denote the set of bilaterally symmetric matchings of $[2n]$.
	Denote by $\mathcal{SCM}_k(n)$ and $\mathcal{SNM}_k(n)$  the sets of $k$-noncrossing and $k$-nonnesting bilaterally symmetric matchings of $\mathcal{SM}(n)$, respectively.

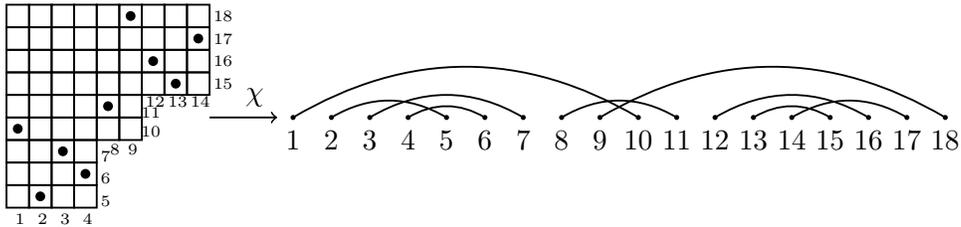
\begin{figure}[H]
		\begin{center}
			\begin{tikzpicture}[font =\small , scale = 0.3,line width = 0.7pt]
				\foreach \i / \j in {1/9,2/9,3/9,4/9,5/6,6/6,7/4,8/4,9/4}
				{
					\foreach \k in {1,...,\j} \draw (\i,-\k)rectangle(\i + 1,-\k -1);
				}
				\foreach \i / \j in {1/6,2/9,3/7,4/8,5/5,6/1,7/3,8/4,9/2} \filldraw[black](\i+0.5,-\j-0.5)circle(5pt);
				
				\draw[->](10,-6)--(13,-6) node[right = -3mm,above]{$\chi$};
				\tikzmath{\c = 1.7;\r = 12;};
				\foreach \x in {1,...,18} \filldraw[black](\x*\c + \r,-6)circle(2pt);
				\foreach \x in {1,...,18} \node at(\x*\c + \r,-7) {\x};
				
				\foreach \x /\y in {1/10,2/5,3/7,4/6,8/11,9/18,12/16,13/15,14/17}
				{
					\draw (\x*\c + \r,-6)..controls(2*\x/3*\c+\y/3*\c+\r,-6+\y/3-\x/3)and
					(\x/3*\c+2*\y/3*\c+\r,-6+\y/3-\x/3)..(\y*\c + \r,-6);
				}			
			    
			    \tikzmath{\r = 0.6;}
				\node[font = \tiny] at (1+\r,-10.5) {$1$};
				\node[font = \tiny] at (2+\r,-10.5) {$2$};
				\node[font = \tiny] at (3+\r,-10.5) {$3$};
				\node[font = \tiny] at (4+\r,-10.5) {$4$};
				\node[font = \tiny] at (5-0.2+\r,-9.7) {$5$};
				\node[font = \tiny] at (5-0.2+\r,-8.7) {$6$};
				\node[font = \tiny] at (5-0.2+\r,-7.7) {$7$};
				\node[font = \tiny] at (5.2+\r,-7.4) {$8$};
				\node[font = \tiny] at (6+\r,-7.4) {$9$};
				\node[font = \tiny] at (5.8+\r+1,-6.6) {$10$};
				\node[font = \tiny] at (5.8+\r+1,-5.8) {$11$};
				\node[font = \tiny] at (6+\r+1,-5.3) {$12$};
				\node[font = \tiny] at (7+\r+1,-5.3) {$13$};
				\node[font = \tiny] at (8+\r+1,-5.3) {$14$};
				\node[font = \tiny] at (9+\r+1,-4.5) {$15$};
				\node[font = \tiny] at (9+\r+1,-3.5) {$16$};
				\node[font = \tiny] at (9+\r+1,-2.5) {$17$};
				\node[font = \tiny] at (9+\r+1,-1.5) {$18$};
			\end{tikzpicture}
		\end{center}
		\caption{ An example of the bijection $\chi$ between $\mathcal{T}_n$ and $\mathcal{M}(n)$.}\label{fig:chi}
	\end{figure}

	Given a matching $M \in \mathcal{M}(n)$ which contains three arcs
	$(i-1,a), (i,b),(i+1,c)$, the index $i$  is said to be a {\em peak} of $M$ if $a<b>c$, whereas $i$ is said to be a {\em valley} of $M$ if $a>b<c$.
	Denote by $\mathrm{Peak}(M)$ and $\mathrm{Val}(M)$ the set of 
	peaks and the set of valleys of $M$, respectively.
	For example, let $M$ be the matching whose linear representation is shown in Figure \ref{fig:chi} (right).
	Then we have $\mathrm{Peak}(M) = \{3\}$ and $\mathrm{Val}(M) = \{2,13\}$.
	The {\em type} of $M$, denoted by $\mathrm{type}(M)$, is the sequence obtained from  $M$ by  tracing from $1$ to $2n$ and writing  $u$ (resp. $d$) whenever we encounter an opener (resp. a closer).
	
	Let $T$ be a transversal in $\mathcal{T}_{\lambda}$.
	The {\em type} of $T$, denoted by $\mathrm{type}(T)$,
	is the sequence obtained from  $T$ by  tracing the south-east  border of $\lambda$ from  south-west to north-east and writing  $u$ (resp. $d$) whenever we encounter a east step (resp. a north step).
	Notice that the type of a transversal $T$ is uniquely determined by its 
	shape $\lambda$.
	
	In the following, we give a description of the bijection $\chi$ from transversals to  matchings  constructed  in \cite{Yan2023}.
	
	\noindent{\bf The bijection $\chi$ from $\mathcal{T}_n$ to $\mathcal{M}(n)$.}\\
	Given a transversal $T \in \mathcal{M}(n)$,
	$M = \chi(T)$ is defined to be the unique matching satisfying that:
	\begin{itemize}
		\item  $\mathrm{type}(M)=\mathrm{type}(T)$;
		\item  There is an arc connecting the $i$-th left-to-right opener and the $j$-th right-to-left closer
		if and only if the square $(i, j)$ is filled with a $1$.
	\end{itemize}
	
	For example, let $T \in \mathcal{T}_{9}$ be the transversal of the Young diagram $\lambda=(9,9,9,9,6,6,4,$ $4,4)$ as shown in Figure \ref{fig:chi} (left).
	By applying  the map $\chi$,  we get a matching $M \in \mathcal{M}(9)$ as shown in Figure \ref{fig:chi} (right).
	It is easily checked that both the types of $T$ and $M$ are $ uuuuddduudduuudddd$.

	Assume that $T$ is a transversal in $\mathcal{T}_{\lambda}$ with $\mathrm{Peak}(T) = \{p_1, p_2, \ldots, p_k\}$
	and $\mathrm{Val}(T) = \{v_1, v_2, \ldots, v_{\ell}\}$.
	We assign labels to the steps in the south-east border of $\lambda$
	with $1,2,\ldots,2n$ from south-west to north-east. 
	If the border corresponding to the bottom of column $p_i$ (resp. $v_i$) receives the label $p_i'$ (resp. $v_i'$), 
	define $\widetilde{\mathrm{Peak}}(T)=\{p_1', p_2', \ldots, p_k'\}$
	and $\widetilde{\mathrm{Val}}(T)=\{v_1', v_2', \ldots, v_{\ell}'\}$.
	For example, for the transversal $T$ in Figure \ref{fig:chi} (left), we have 
	$\mathrm{Peak}(T) = \{2,8\}$, $\widetilde{\mathrm{Peak}}(T) = \{2,13\}$, $\mathrm{Val}(T) = \{3\}$ and $\widetilde{\mathrm{Val}}(T) = \{3\}$.
	It is apparent that for a transversal $T$ of a given shape,
	$\widetilde{\mathrm{Peak}}(T)$ (resp. $\widetilde{\mathrm{Val}}(T)$) is uniquely determined by $\mathrm{Peak}(T)$ (resp. $\mathrm{Val}(T)$), and vise versa.
	Then the following theorem can be easily summarized from the relevant work in \cite{Yan2023,Zhou}.
	
	\begin{theorem}(\cite{Yan2023,Zhou})\label{thm:chi}
		The map $\chi$ is a bijection between $\mathcal{T}_n$ and $\mathcal{M}(n)$
		such that  for any  $T\in \mathcal{T}_n$, its corresponding matching  $M = \chi(T)$ verifies that
		\begin{enumerate}[label=\upshape(\roman*)]
			\item  $\mathrm{type}(T)=\mathrm{type}(M)$;	
			\item  $\widetilde{\mathrm{Peak}}(T) = \mathrm{Val}(M)$;
			\item  $\widetilde{\mathrm{Val}}(T) = \mathrm{Peak}(M)$;
			\item  $T\in \mathcal{ST}_n$ if and only if $M \in \mathcal{SM}(n)$;
			\item  $T\in \mathcal{T}_n(J_k)$ if and only if $M\in \mathcal{CM}_k(n)$;
			\item  $T\in \mathcal{T}_n(I_k)$ if and only if $M\in \mathcal{NM}_k(n)$.
		\end{enumerate}
	\end{theorem}
	
	\subsection{The bijection $\phi$ from matchings to oscillating tableaux}
	An {\em oscillating    tableau} of shape $\lambda$ and length $n$ is a sequence $(\emptyset=\lambda^0, \lambda^1, \ldots, \lambda^{n}=\lambda)$ of integer partitions such that     $\lambda^{i}$ is obtained from $\lambda^{i-1}$ by  either adding a square or deleting a square.
	In what follows,  oscillating  tableaux are always of shape $\emptyset$ unless specified otherwise. 
	
	In \cite{Sundaram}, Sundaram constructed  a bijection $\phi$ between matchings  and  oscillating  tableaux. 
	The explicit description of the bijection between matchings and  oscillating    tableaux can also be found in \cite{ChenD} and \cite{StanleyVol2}.
	Chen-Deng-Du-Stanley-Yan \cite{ChenD} proved that the bijection $\phi$ verifies  the following celebrated property.
	\begin{lemma}\label{lem:chenth6}
		(\cite{ChenD},  Theorem 6)
		Let $M$ be a matching of $[2n]$ and $\phi(M)=(\lambda^0, \lambda^1, \ldots, $ $\lambda^{2n})$.  
		Then	$M$ is $k$-noncrossing (resp. $k$-nonnesting) if and only if  any $\lambda^i$ contains at most $k-1$
		rows (resp. columns).
	\end{lemma}

	For a sequence of integer partitions $U = (\mu^1, \mu^2,\ldots, \mu^k)$,
	define $U^r = (\mu^k,\ldots, \mu^2, \mu^1)$.
	We call $U^r$ the {\em reverse} of $U$.
	Given an oscillating tableau $O$,  $O^r$ is still an oscillating tableau.
	If $O^{r}=O$, we say that $O$ is {\em symmetric}.
	By using Sch\"{u}tzenberger's theorem (see \cite{StanleyVol2}, Chapter 7.11) for the RSK correspondence (see \cite{StanleyVol2}, Chapter 7.13), Xin-Zhang \cite{xinzhang} proved that the bijection $\phi$ has the following celebrated property.
	
	\begin{lemma}\label{lem:xinzhang}
		(\cite{xinzhang},  Theorem 1)
		For any given matching $M$ and oscillating tableau $O$, $\phi(M^{r})=O^{r}$ if and only if $\phi(M)=O$.
	\end{lemma}
	
	Let $\mathcal{O}(n)$ and $\mathcal{SO}(n)$ denote the set of 
	oscillating tableaux and symmetric oscillating tableaux of length $2n$,
	respectively. 
	Let $\mathcal{OR}_{k}(n)$ (resp. $\mathcal{OC}_{k}(n)$) denote the set of  oscillating  tableaux $(\lambda^0, \lambda^1, \ldots, $ $\lambda^{2n})$ of length $2n$ in which any $\lambda^i$ has at most $k-1$ rows (resp. columns).
	Let $\mathcal{SOR}_{k}(n)$ (resp. $\mathcal{SOC}_{k}(n)$) denote the set of  symmetric oscillating  tableaux of $\mathcal{OR}_{k}(n)$ (resp. $\mathcal{OC}_{k}(n)$).

	\noindent {\bf The bijection $\phi$ from matchings to oscillating Tableaux.}\\
	Given a matching $M \in \mathcal{M}(n)$, we will recursively define a sequence of 
	SYT's $T_0,T_1,\ldots,T_{2n}$ as follows: Start from the empty SYT by letting $T_{2n}=\emptyset$, read the number $j\in [2n]$ one by one from $2n$ to 1, and let $T_{j-1}$ be the SYT obtained from $T_j$ for each $j$ by the following procedure.
	\begin{itemize}
		\item If  $j$ is the closer of an arc $(i,j)$, then insert $i$ (by the RSK algorithm) into the tableau $T_j$.
		\item If $j$ is the opener of an arc $(j, k)$,  then  remove $j$ from the tableau $T_j$.
	\end{itemize}
	Then the oscillating tableau $\phi(M)$  is the sequence of shapes of the above SYT's.
	
	For example, let $M$ be the  matching as shown in Figure \ref{fig:chi} (right).
	Table \ref{table:phi} describes in detail the processes when applying $\phi$ to
	$M$ to obtain the corresponding oscillating tableau $O = (\lambda^0, \lambda^1, \ldots, \lambda^{18})$.

	\begin{table}[H]
		\centering
		\scriptsize
		\renewcommand\arraystretch{1}
		\caption{An oscillating tableau $\phi(M)$ corresponding to the matching  $M=\{(1, 10), (2,5), (3,7), (4, 6), (8,11), (9,18), (12,16), (13,15), (14,17)\}$.}\label{table:phi}
		\vskip 2mm
		
		\begin{tabular}{|c|m{2cm}<{\centering}|c|c| m{2cm}<{\centering}|c|}
			\hline
			$i$ &  \vskip 0.5mm  $T_i$       & $\lambda^i$      &  $i$   & \vskip 0.5mm $T_i$     & $\lambda^i$    \\ \hline
			$0$ &  $\emptyset$   & $\emptyset$       & $10$    &
			\vskip 1mm
			$\begin{array}[b]{*{1}c}\cline{1-1}
				\lr{8}\\\cline{1-1}
				\lr{9}\\\cline{1-1}
			\end{array}$
			& $(1,1)$                                    \\ \hline
			
			$1$  &
			$\begin{array}[b]{*{1}c}\cline{1-1}
				\lr{1}\\\cline{1-1}
			\end{array}$
			& $(1)$ & $11$ &  \vskip 1mm
			$\begin{array}[b]{*{1}c}\cline{1-1}
				\lr{9}\\\cline{1-1}
			\end{array}$
			& $(1)$                                 \\ \hline
			
			$2$  & \vskip 1mm
			$\begin{array}[b]{*{2}c}\cline{1-2}
				\lr{1}&\lr{2}\\\cline{1-2}
			\end{array}$
			& $(2)$  & $12$ &  \vskip 1mm
			$\begin{array}[b]{*{2}c}\cline{1-2}
				\lr{9}&\lr{12}\\\cline{1-2}
			\end{array}$
			&   $(2)$                             \\ \hline
			
			$3$  & \vskip 1mm
			$\begin{array}[b]{*{2}c}\cline{1-2}
				\lr{1}&\lr{2}\\\cline{1-2}
				\lr{3}\\\cline{1-1}
			\end{array}$
			& $(2,1)$  & $13$  & \vskip 1mm
			$\begin{array}[b]{*{3}c}\cline{1-3}
				\lr{9}&\lr{12}&\lr{13}\\\cline{1-3}
			\end{array}$
			&    $(3)$                            \\ \hline
			
			$4$  & \vskip 1mm
			$\begin{array}[b]{*{3}c}\cline{1-3}
				\lr{1}&\lr{2}&\lr{4}\\\cline{1-3}
				\lr{3}\\\cline{1-1}
			\end{array}$
			& $(3,1)$  & $14$  & \vskip 1mm
			$\begin{array}[b]{*{3}c}\cline{1-3}
				\lr{9}&\lr{12}&\lr{13}\\\cline{1-3}
				\lr{14}\\\cline{1-1}
			\end{array}$
			&    $(3,1)$                            \\ \hline
			
			$5$  & \vskip 1mm
			$\begin{array}[b]{*{3}c}\cline{1-3}
				\lr{1}&\lr{3}&\lr{4}\\\cline{1-3}
			\end{array}$
			& $(3)$  & $15$  & \vskip 1mm
			$\begin{array}[b]{*{2}c}\cline{1-2}
				\lr{9}&\lr{12}\\\cline{1-2}
				\lr{14}\\\cline{1-1}
			\end{array}$
			&    $(2,1)$                            \\ \hline
			
			$6$  & \vskip 1mm
			$\begin{array}[b]{*{2}c}\cline{1-2}
				\lr{1}&\lr{3}\\\cline{1-2}
			\end{array}$
			& $(2)$  & $16$  & \vskip 1mm
			$\begin{array}[b]{*{2}c}\cline{1-2}
				\lr{9}&\lr{14}\\\cline{1-2}
			\end{array}$
			&    $(2)$                            \\ \hline
			
			$7$  & \vskip 1mm
			$\begin{array}[b]{*{1}c}\cline{1-1}
				\lr{1}\\\cline{1-1}
			\end{array}$
			& $(1)$  & $17$  & \vskip 1mm
			$\begin{array}[b]{*{1}c}\cline{1-1}
				\lr{9}\\\cline{1-1}
			\end{array}$
			&  $(1)$        \\ \hline
			
			$8$  & \vskip 1mm
			$\begin{array}[b]{*{2}c}\cline{1-1}
				\lr{1}\\\cline{1-1}
				\lr{8}\\\cline{1-1}
			\end{array}$
			& $(1,1)$  & $18$  & $\emptyset$  &  $\emptyset$       \\ \hline
			
			$9$  & \vskip 1mm
			$\begin{array}[b]{*{3}c}\cline{1-1}
				\lr{1}\\\cline{1-1}
				\lr{8}\\\cline{1-1}
				\lr{9}\\\cline{1-1}
			\end{array}$
			& $(1,1,1)$  &   &   &         \\ \hline
		\end{tabular}
	\end{table}

	Given an oscillating tableau $O=(\lambda^0, \lambda^1, \ldots, \lambda^n)$,  we can associate  it with  a sequence  of $u$'s and $d$'s obtained from  $O$ by   reading $O$  forward and   writing   $u$  (resp.   $d$) whenever  $\lambda^i$ is obtained from $\lambda^{i-1}$ by adding (resp. deleting) a square. Such a sequence is called the {\em type} of $O$, denoted by $\mathrm{type}(O)$.
	Let $O=(\lambda^0, \lambda^1, \ldots, \lambda^n)$ be the oscillating tableau as shown in Table \ref{table:phi}.
	Then we have 
	$\mathrm{type}(O) = uuuuddduudduuudddd$.

	For an oscillating tableau
	$O=(\lambda^0, \lambda^1, \ldots, \lambda^{2n})\in \mathcal{O}(n)$, 
	we associate $O$ with a word $y = y_1y_2\cdots y_{2n} = \bar{\alpha}(O)$ of length $2n$ where $y_i = i$ (resp. $\bar{i}$) if $\lambda^i$ is obtained from 
	$\lambda^{i-1}$ by adding (resp. removing) one square in row $i$.
	We also call $y$ the Yamanouchi word of $O$ as 
	we will see later it is closely related to the Yamanouchi words of SYT's.
	
	Let $O$ be an oscillating tableau of length $2n$ and let $y=y_1y_2\cdots y_{2n}$ be the corresponding Yamanouchi word of $O$.
	An index  $i$ ($1 < i < 2n$) 
	is said to be a {\em peak } of $O$ if  $y_{i-1},y_i,y_{i+1}$ have no bars above them and  $y_{i-1} < y_i \geq y_{i+1}$. Similarly,   an index $i$ ($1 < i < 2n$) is said to be a {\em valley} of $O$ if  $y_{i-1},y_i,y_{i+1}$ have no bars above them  and $y_{i-1} \geq  y_i < y_{i+1}$. 
	Let $\mathrm{Peak}(O)$ and $\mathrm{Val}(O)$ denote the set of peaks and the set of valleys of $O$, respectively.
	For instance, let $O = (\emptyset, (1),(1,1),(2,1),(2,2),(3,2),(3,1),(2,1),$ $(2),(1),\emptyset)$ be an oscillating tableau of length $10$.
	Then we have $y = \bar{\alpha}(O) = 12121\bar{2}\bar{1}\bar{2}\bar{1}\bar{1}$, $\mathrm{Val}(O) = \{3\}$  and $\mathrm{Peak}(O) = \{2,4\}$.
	
	To establish a correspondence between the peak set (or valley set) of 
	matchings and the peak set (or valley set) of oscillating tableaux,
	we shall introduce an equivalent description of $\phi$.
	
	\noindent {\bf The Bijection $\bar{\phi}$ from Matchings to Oscillating Tableaux.}\\
	Given a matching $M \in \mathcal{M}(n)$, we will recursively define a sequence of 
	SYT's $T_0,T_1,\ldots,T_{2n}$ as follows: Start from the empty SYT by letting $T_{0}=\emptyset$, read the number $j\in [2n]$ one by one from 1 to $2n$, and let $T_{j}$ be the SYT obtained from $T_{j-1}$ for each $j$ by the following procedure.
	\begin{itemize}
		\item If  $j$ is the opener of an arc $(j,k)$, then insert $2n+1-k$ (by the RSK algorithm) into the tableau $T_{j-1}$.
		\item If $j$ is the closer of an arc $(i,j)$,  then remove $2n+1-j$ from the tableau $T_{j-1}$.
	\end{itemize}
	Then the oscillating tableau $\bar{\phi}(M)$  is the sequence of shapes of the above SYT's.
	One can easily check that $\bar{\phi}(M) = (\phi(M^r))^r$.
	Hence $\bar{\phi}$ is well-defined.
	Moreover,  by Lemma \ref{lem:xinzhang}, we derive that 
	$\bar{\phi}(M) = (\phi(M^r))^r = (\phi(M)^r)^r = \phi(M)$.
	In what follows, we will treat $\phi$ and $\bar{\phi}$ as identical.
	Now we proceed to    show that the bijection $\phi$ has the following desired properties.

	\begin{theorem}\label{thm:phi}
		The map $\phi$ is a bijection between $\mathcal{M}(n)$ and $\mathcal{O}(n)$
		such that  for any  $M\in \mathcal{M}(n)$, its corresponding oscillating tableau  $O = \phi(M)$ verifies that
		\begin{enumerate}[label=\upshape(\roman*)]
			\item  $\mathrm{type}(M)=\mathrm{type}(O)$;	
			\item  $\mathrm{Peak}(M) = \mathrm{Peak}(O)$;
			\item $\mathrm{Val}(M) = \mathrm{Val}(O)$;
			\item  $M\in \mathcal{SM}(n)$ if and only if $O \in \mathcal{SO}(n)$;
			\item  $M\in \mathcal{CM}_k(n)$ if and only if $O\in \mathcal{OR}_k(n)$;
			\item  $M\in \mathcal{NM}_k(n)$ if and only if $O\in \mathcal{OC}_k(n)$.
		\end{enumerate}
	\end{theorem}
	
	\pf
	(\upshape \rmnum{1}), (\upshape \rmnum{4}), (\upshape \rmnum{5}) and (\upshape \rmnum{6}) follow directly from  Lemmas \ref{lem:chenth6} and \ref{lem:xinzhang}. 
	Now we will only deal with (\upshape \rmnum{2}), and 
	(\upshape \rmnum{3}) can be deduced by similar arguments. 
	
	If $i \in \mathrm{Peak}(M)$, then there exist three arcs $(i-1,a),(i,b),(i+1,c)$ of 
	$M$ such that $a<b>c$.
	Let $y=y_1y_2\cdots y_{2n}$ be the corresponding Yamanouchi word of $O$.
	By the definition of $\bar{\phi}$ and $\bar{\alpha}$, we have that 
	$y_{i-1},y_i,y_{i+1}$ record the added row indices when $2n+1-a,2n+1-b,2n+1-c$ 
	are inserted one by one to a tableau.
	Note that $a<b>c$, we have $2n+1-a>2n+1-b<2n+1-c$.
	By the property of RSK algorithm, 
	we deduce that $y_{i-1}<y_i\geq y_{i+1}$, namely, $i \in \mathrm{Peak}(O)$.
	Thus $\mathrm{Peak}(M) \subseteq \mathrm{Peak}(O)$.
	Since the above process is reversible, we have $\mathrm{Peak}(O) \subseteq \mathrm{Peak}(M)$.
	Hence $\mathrm{Peak}(M) = \mathrm{Peak}(O)$, this completes the proof of (\upshape \rmnum{2}).
	\qed

	\subsection{An involution $\Phi$ on transversals}
	In this subsection, we aim to establish an involution $\Phi$ on $\mathcal{T}_n$.  To this end, we shall
	construct an involution $\gamma$ on $\mathcal{O}(n)$ .
	
%

	Given an oscillating tableau $O=(\lambda^0, \lambda^1, \ldots, $ $\lambda^{2n})$,
	define $\gamma(O) = ((\lambda^0)^T, (\lambda^1)^T, \ldots,$ $(\lambda^{2n})^T)$.
	We call $\gamma(O)$ the {\em conjugate} of $O$.
	It is easily seen that $\gamma(O)$ is also an oscillating tableau.

	\begin{lemma}\label{lem:conjugate}
		The conjugate map $\gamma$  is an involution on $\mathcal{O}(n)$ 
		such that  for any  $O\in \mathcal{O}(n)$, its corresponding oscillating tableau  $O' = \gamma(O)$ verifies that 
		\begin{enumerate}[label=\upshape(\roman*)]
			\item  $\mathrm{type}(O)=\mathrm{type}(O')$;	
			\item  $\mathrm{Val}(O) = \mathrm{Peak}(O')$;
			\item  $\mathrm{Peak}(O) = \mathrm{Val}(O')$;
			\item  $O\in \mathcal{SO}(n)$ if and only if $O' \in \mathcal{SO}(n)$;
			\item  $O\in \mathcal{OR}_k(n)$ if and only if $O'\in \mathcal{OC}_k(n)$;
			\item  $O\in \mathcal{OC}_k(n)$ if and only if $O'\in \mathcal{OR}_k(n)$.
		\end{enumerate}
	\end{lemma}
	
	\pf
	(\upshape \rmnum{1}),  (\upshape \rmnum{4}),
	(\upshape \rmnum{5}) and (\upshape \rmnum{6}) follow directly from the definition of the conjugate of an oscillating tableau $O$.

	Given an oscillating tableau $O=(\lambda^0, \lambda^1, \ldots, \lambda^{2n})$
	with its corresponding Yamanouchi word $y = y_1y_2\cdots y_{2n}$.
	Let $x = x_1x_2\cdots x_{2n}$ be the word of length $2n$ where
	$x_i = i$ (resp. $\bar{i}$) if $\lambda^i$ is obtained from 
	$\lambda^{i-1}$ by adding (resp. removing) one square in column $i$.
	Then by the definition of the conjugate of an oscillating tableau $O$, 
	we deduce that  $x$ is indeed the  Yamanouchi word of $O'$.
	Let $y_i,y_{i+1}$ be two adjacent numbers of $y$ without bars above them.
	Then $x_i,x_{i+1}$ also have no bars above them.
	It is easily seen  that 
	$y_i < y_{i+1}$ if and only if $x_i \geq x_{i+1}$, and 
	$y_i \geq y_{i+1}$ if and only if $x_i < x_{i+1}$.
	Then the (\upshape \rmnum{2}) and  (\upshape \rmnum{3}) follow directly from 
	the definitions of peak and valley of oscillating tableaux.
	\qed

	Notice that for any two transversals $T_1,T_2$  of the same shape  $\lambda$
	and any  $\mathrm{stat}_1,\mathrm{stat}_2 \in \{\mathrm{Peak},\mathrm{Val}\}$,
	we have $\widetilde{\mathrm{stat}}_1(T_1) = \widetilde{\mathrm{stat}}_2(T_2)$ if and only if $\mathrm{stat}_1(T_1) = \mathrm{stat}_2(T_2)$.
	Recall that the type of a transversal is uniquely determined by its shape. By Theorem \ref{thm:chi},  Theorem  \ref{thm:phi} and Lemma \ref{lem:conjugate},  we deduce the following result.
	\begin{theorem}\label{lem:Phi}
		The map $\Phi=\chi^{-1}\circ \phi^{-1}\circ\gamma\circ    \phi \circ \chi$ induces an involution on $\mathcal{T}_n$
		such that for any $T\in \mathcal{T}_n$,
		its corresponding transversal $S = \Phi(T)$ verifies that
		\begin{enumerate}[label=\upshape(\roman*)]
			\item  $\mathrm{type}(T)=\mathrm{type}(S)$;	
			\item  $\mathrm{Peak}(T) = \mathrm{Val}(S)$;
			\item  $\mathrm{Val}(T) = \mathrm{Peak}(S)$;
			\item  $T\in \mathcal{ST}_n$ if and only if $S \in \mathcal{ST}_n$;
			\item  $T\in \mathcal{T}_n(J_k)$ if and only if $S\in \mathcal{T}_n(I_k)$;
			\item  $T\in \mathcal{T}_n(I_k)$ if and only if $S\in \mathcal{T}_n(J_k)$.
		\end{enumerate}
	\end{theorem}

	When restricted to transversals in $\mathcal{T}_{\lambda}$ for any given Young diagram $\lambda$,
	$\Phi$ induces an involution on $\mathcal{T}_{\lambda}$. 
	Then we have the following corollary of Theorem  \ref{lem:Phi}.
	
	\begin{corollary}\label{coro:T}
		For any Young diagram $\lambda$,  $\mathrm{Peak}$ and $\mathrm{Val}$ have a
		symmetric joint distribution over $\mathcal{T}_{\lambda}$, that is, 
		$$\sum_{T \in \mathcal{T}_{\lambda} }t^{\mathrm{Peak}(T)}q^{\mathrm{Val}(T)}
		=\sum_{T \in \mathcal{T}_{\lambda}}t^{\mathrm{Val}(T)}q^{\mathrm{Peak}(T)}.$$
	\end{corollary}
	
	When further restricted to transversals in $\mathcal{ST}_{\lambda}$ for any self-conjugate Young diagram $\lambda$,
	$\Phi$ induces an involution on $\mathcal{ST}_{\lambda}$.
	Then we have the following corollary of Theorem \ref{lem:Phi}.
	
	\begin{corollary}\label{coro:ST}
		For any self-conjugate Young diagram $\lambda$, $\mathrm{Peak}$ and $\mathrm{Val}$ have a 
		symmetric joint distribution over $\mathcal{ST}_{\lambda}$, that is, 
		$$\sum_{T \in \mathcal{ST}_{\lambda} }t^{\mathrm{Peak}(T)}q^{\mathrm{Val}(T)}
		=\sum_{T \in \mathcal{ST}_{\lambda}}t^{\mathrm{Val}(T)}q^{\mathrm{Peak}(T)}.$$
	\end{corollary}
	
	When  restricted to transversals in $\mathcal{T}_{\lambda}(J_k)$ for any given Young diagram $\lambda$ and for any positive integer $k$, 
	we have the following refinements of Corollaries \ref{coro:T} and \ref{coro:ST}, respectively.
	
	\begin{corollary}\label{coro:TJk}
		For any Young diagram $\lambda$ and any positive integer $k$,  we have
		$$\sum_{T \in \mathcal{T}_{\lambda}(J_k) }t^{\mathrm{Peak}(T)}q^{\mathrm{Val}(T)}
		=\sum_{T \in \mathcal{T}_{\lambda}(I_k)}t^{\mathrm{Val}(T)}q^{\mathrm{Peak}(T)}.$$
	\end{corollary}
	
	\begin{corollary}\label{coro:STJk}
		For any self-conjugate Young diagram $\lambda$ and any positive integer $k$, we have 
		$$\sum_{T \in \mathcal{ST}_{\lambda}(J_k)}t^{\mathrm{Peak}(T)}q^{\mathrm{Val}(T)}
		=\sum_{T \in \mathcal{ST}_{\lambda}(I_k)}t^{\mathrm{Val}(T)}q^{\mathrm{Peak}(T)}.$$
	\end{corollary}

	Setting  $q_i = 1$ for $i\geq 1$  in Corollaries \ref{coro:TJk} and \ref{coro:STJk}, we derive that
	\begin{equation}\label{equ:1}
		\sum_{T \in \mathcal{T}_{\lambda}(J_k) }t^{\mathrm{Peak}(T)}
		=\sum_{T \in \mathcal{T}_{\lambda}(I_k)}t^{\mathrm{Val}(T)}
	\end{equation}
	and \begin{equation}\label{equ:3}
		\sum_{T \in \mathcal{ST}_{\lambda}(J_k) }t^{\mathrm{Peak}(T)}
		=\sum_{T \in \mathcal{ST}_{\lambda}(I_k)}t^{\mathrm{Val}(T)}.
	\end{equation}
	Setting  $t_i = 1$ for $i\geq 1$  in  Corollaries \ref{coro:TJk} and \ref{coro:STJk}, we derive that
	\begin{equation}\label{equ:2}
		\sum_{T \in \mathcal{T}_{\lambda}(J_k) }q^{\mathrm{Val}(T)}
		=\sum_{T \in \mathcal{T}_{\lambda}(I_k)}q^{\mathrm{Peak}(T)}
	\end{equation}
	and
	\begin{equation}\label{equ:4}
		\sum_{T \in \mathcal{ST}_{\lambda}(J_k) }q^{\mathrm{Val}(T)}
		=\sum_{T \in \mathcal{ST}_{\lambda}(I_k)}q^{\mathrm{Peak}(T)}.
	\end{equation}

	\subsection{A bijection $\Psi$ from $\mathcal{T}_n$ to itself }
	In this subsection, we aim to establish a bijection  $\Psi$ from $\mathcal{T}_n$ to itself which enables us to prove Theorem \ref{thm:transversal} and Theorem \ref{thm:symm-transversal}. 
	To this end, we shall construct a bijection $\psi$ from $\mathcal{O}(n)$ to itself relying on Theorem \ref{thm:skewSYT}.

	For an oscillating tableau  $O=(\lambda^0, \lambda^1, \ldots, \lambda^{2n})$, 
it can be uniquely decomposed  as $(A_1,D_1,A_2,D_2,\ldots ,A_k,D_k)$,
where  each $A_i$  (resp. $D_i$) is called an {\em addition run}  (resp. a {\em deletion run}) of $O$ consisting of  a maximal chain of consecutive  partitions such that each partition is obtained from 
the one before it by adding (resp. removing) one square. 	
	 In what follows, such a decomposition is called the {\em addition-deletion decomposition} of $O$.   
	For example, let $O=(\lambda^0, \lambda^1, \ldots, \lambda^{18})$ be the oscillating tableau as shown in Table \ref{table:phi}.
	Then the addition-deletion decomposition of $O$ is given by $((\lambda^0,\lambda^1,\lambda^2,\lambda^3,\lambda^4), (\lambda^4,\lambda^5,\lambda^6,\lambda^7),(\lambda^7,\lambda^8,\lambda^9), (\lambda^9,\lambda^{10},\lambda^{11}), (\lambda^{11},\lambda^{12},\lambda^{13},\lambda^{14}), (\lambda^{14},\lambda^{15}, $   $ \lambda^{16},\lambda^{17},\lambda^{18}))$.

	\begin{observation} \label{obs1}
		Given an   oscillating tableau  $(\lambda^0, \lambda^1, \ldots, \lambda^{2n})$, assume that the addition-deletion decomposition of  $O$ is given by   $(A_1,D_1,A_2,D_2,\ldots ,A_k,D_k)$. 
		Then $O \in \mathcal{SO}(n)$ if and only if
		$A_i = D_{k+1-i}^r$ for all $1\leq i\leq k$.   
	\end{observation}
	
	Given a skew diagram $\lambda/\mu$ of size $n$, let
	$\mathcal{AR}(\lambda/\mu)$ denote the set of  sequences $A=(\mu=\lambda^{0},\lambda^{1},\ldots, \lambda^{n}=\lambda)$ of partitions in which   $\lambda^{i}$ is obtained from $\lambda^{i-1}$ by adding a square. 
	Let $\mathcal{DR}(\lambda/\mu)$ denote the set of  sequences  $D=(\lambda=\lambda^{0},\lambda^{1},\ldots, \lambda^{n}=\mu)$ of partitions in which   $\lambda^{i}$ is obtained from $\lambda^{i-1}$ by deleting a square.

	Now we build a bijection $\theta: \mathcal{AR}(\lambda/\mu) \rightarrow \mathcal{Y}(\lambda/\mu)$. 
	Given a sequence $ A=(\mu=\lambda^{0},\lambda^{1},\ldots,$ $ \lambda^{n}=\lambda)\in  \mathcal{AR}(\lambda/\mu)$, assume that
	$\lambda^{i}$ is obtained from $\lambda^{i-1}$ by adding a square at row $y_i$ for all $1\leq i \leq n$, then set $\theta(A) = y_1y_2\ldots y_n$.
	Clearly,  $\theta(A)\in \mathcal{Y}(\lambda/\mu)$. 
	Conversely, given a Yamanouchi word  $y=y_1y_2\ldots y_n\in \mathcal{Y}_{\lambda/\mu} $, one can easily  recover a sequence $ (\mu=\lambda^{0},\lambda^{1},\ldots, \lambda^{n}=\lambda)$ of partitions, denoted by $\theta'(y)$,  in which $\lambda^i$ is obtained from $\lambda^{i-1}$ by adding a square located  at row $y_i$. 
	It is easily seen that the maps $\theta$ and $\theta'$ are inverses of each other and hence $\theta $ is a bijection. 
	Define $\bar{\theta}(D) = \theta(D^r)$ for any $D\in \mathcal{DR}(\lambda/\mu)$.
	It is easily checked that $\bar{\theta}$ induces a bijection between $\mathcal{DR}(\lambda/\mu)$ and 
	$\mathcal{Y}(\lambda/\mu)$.

	\begin{lemma} \label{lem:psi}
		Let $n\geq 1$. There exists a bijection $\psi: \mathcal{O}(n) \rightarrow \mathcal{O}(n)$ such that for any oscillating tableau $O\in \mathcal{O}(n)$, its corresponding oscillating tableau  $O'=\psi(O)$ verifies that
		\begin{enumerate}[label=\upshape(\roman*)]
			\item  $\mathrm{type}(O)=\mathrm{type}(O')$;	
			\item  $\mathrm{Val}(O) = \mathrm{Peak}(O')$;
			\item $O \in \mathcal{SO}(n)$ if and only $O' \in \mathcal{SO}(n)$;
			\item  $O\in \mathcal{OR}_k(n)$ if and only if $O'\in \mathcal{OR}_k(n)$;
			\item  $O\in \mathcal{OC}_k(n)$ if and only if $O'\in \mathcal{OC}_k(n)$.
		\end{enumerate}

	\end{lemma}
	\pf 
	Now we first give a description of   the map $\psi$.
	Given an   oscillating tableau  $O=(\lambda^0, \lambda^1, \ldots, \lambda^{2n})$, assume that the addition-deletion decomposition of  $O$ is given by   $$(A_1,D_1,A_2,D_2,\ldots ,A_k,D_k).$$ 
	Let $W_i = \theta(A_i)$ and $V_i = \bar{\theta}(D_i)$ for all $1\leq i\leq k$.
	By Theorem  \ref{thm:y}, for any skew shape $\lambda/ \mu$,  there exists a bijection $  \xi: \mathcal{Y}(\lambda/ \mu)\rightarrow   \mathcal{Y}(\lambda/ \mu)$ such that for any $w\in \mathcal{Y}(\lambda/ \mu)\ $, we have $\mathrm{Val}(w) = \mathrm{Peak}(\xi(w))$. 
	Let $W'_i=\xi(W_i)$ and $V'_i=\xi(V_i)$ for all $1\leq i\leq k$. 
	Define  $O'=\psi(O)$  to be the  oscillating tableau whose addition-deletion decomposition is given by  $$(A'_1,D'_1,A'_2,D'_2,\ldots ,A'_k, D'_k),$$ where    $A'_i= \theta^{-1}(W'_i)$  and $D'_i= (\bar{\theta}^{-1}(V'_i))$ for $1\leq i\leq k$. 
	
	By the definitions of $\theta$, $\bar{\theta}$ and 
	$\xi$, one can easily find that 
	$A_i\in \mathcal{AR}(\lambda/\mu)$ if and only if $A'_i\in \mathcal{AR}(\lambda/\mu)$ and 
	$D_i\in \mathcal{DR}(\lambda/\mu)$ if and only if $D'_i\in \mathcal{DR}(\lambda/\mu)$ for all $1\leq i\leq k$.
	This implies that $O'$ is indeed an oscillating tableau.
	Moreover,  by the construction of $\psi$,  the map $\psi$ verifies the properties  (\upshape \rmnum{1}), (\upshape \rmnum{4}) and (\upshape \rmnum{5}). 
	Notice that 
	$\mathrm{Val}(O)$ is uniquely determined by the sets $\mathrm{Val}(W_i)$ for all $1\leq i\leq k$ and $\mathrm{Peak}(O')$ is uniquely determined by the sets $\mathrm{Peak}(W'_i)$ for all $1\leq i\leq k$.  Since  $\mathrm{Val}(W_i)=\mathrm{Peak}(W'_i)$ for all $1\leq i\leq k$. This yields that $\mathrm{Val}(O)=\mathrm{Peak}(O')$ as desired, completing the proof of (\upshape \rmnum{2}). 
	
	It is routine to check that
	\begin{align*}
		O \in \mathcal{SO}(n) &\Leftrightarrow A_i = D_{k+1-i}^r\\ &\Leftrightarrow  W_i =\theta(A_i) = \theta(D_{k+1-i}^r) =V_{k+1-i} \\
		&\Leftrightarrow W'_i  = \xi(W_i)= \xi (V_{k+1-i}) = V'_{k+1-i} \\
		&\Leftrightarrow A'_i = \theta^{-1} (W'_i)=\theta^{-1}(V'_{k+1-i})= (D'_{k+1-i})^r\\
		&\Leftrightarrow O'\in \mathcal{SO}(n)
	\end{align*}
	for all $1\leq i\leq k$. Thus  (\upshape \rmnum{3}) follows, completing the proof. 
	\qed

	By Theorem \ref{thm:chi},  Theorem  \ref{thm:phi} and Lemma \ref{lem:psi},  we deduce the following result.
	\begin{theorem}\label{thm:Psi}
		The map $\Psi=\chi^{-1}\circ \phi^{-1} \circ \psi  \circ \phi \circ \chi$  induces a bijection from   $\mathcal{T}_n$ to itself
		such that for any  $T\in \mathcal{T}_n$, its corresponding transversal  $S = \Psi(T)$ verifies that
		\begin{enumerate}[label=\upshape(\roman*)]
			\item  $\mathrm{type}(T)=\mathrm{type}(S)$;	
			\item  $\mathrm{Peak}(T) = \mathrm{Val}(S)$;
			\item $ T \in \mathcal{ST}_n$ if and only if 
			$S \in \mathcal{ST}_n$;
			\item  $T\in \mathcal{T}_n(J_k)$ if and only if $S\in \mathcal{T}_n(J_k)$;
			\item  $T\in \mathcal{T}_n(I_k)$ if and only if $S\in \mathcal{T}_n(I_k)$.
		\end{enumerate}
	\end{theorem}
	
\noindent  {\bf Proof of Theorem \ref{thm:transversal}.} 
From (\upshape \rmnum{1}), (\upshape \rmnum{2}) and (\upshape \rmnum{4}) in Theorem
\ref{thm:Psi}, 
we deduce that  Peak and Val are 
equidistributed over the set $\mathcal{T}_{\lambda}(J_k)$, namely,
$$\sum_{T \in \mathcal{T}_{\lambda}(J_k) }t^{\mathrm{Val}(T)}
= \sum_{T \in \mathcal{T}_{\lambda}(J_k) }t^{\mathrm{Peak}(T)}.$$
Combining (\ref{equ:1}) and (\ref{equ:2}), Theorem \ref{thm:transversal}
follows.
\qed

\noindent  {\bf Proof of Theorem \ref{thm:symm-transversal}.} 
From (\upshape \rmnum{1}), (\upshape \rmnum{2}), (\upshape \rmnum{3}) and (\upshape \rmnum{4}) in Theorem \ref{thm:Psi}, 
we deduce that  Peak and Val are 
equidistributed over the set $\mathcal{ST}_{\lambda}(J_k)$, namely,
$$\sum_{T \in \mathcal{ST}_{\lambda}(J_k) }t^{\mathrm{Val}(T)}
= \sum_{T \in \mathcal{ST}_{\lambda}(J_k) }t^{\mathrm{Peak}(T)}.$$
Combining (\ref{equ:3}) and (\ref{equ:4}), Theorem \ref{thm:symm-transversal}
follows.
\qed

	\subsection{Proofs of Theorems \ref{thm:transversal-general} and \ref{thm:symm-transversal-general}}
	
	In \cite{BWX}, Backelin-West-Xin obtained the following  lemma.
	\begin{lemma}(\cite{BWX}, Proposition 2.3)\label{lem:BWX}
		For  any  permutations $\alpha$, $\beta$ and $\tau$, if $|\mathcal{T}_{\mu}(\alpha)| =  |\mathcal{T}_{\mu}(\beta)|$ for any Young diagram $\mu$,  then  $|\mathcal{T}_{\lambda}(\alpha\oplus \tau)| =  |\mathcal{T}_{\lambda}(\beta\oplus \tau)|$ for any Young diagram $\lambda$.
	\end{lemma}

	In the spirit of Backelin-West-Xin's  proof of Lemma \ref{lem:BWX}, we  deduce the following 
	theorem, which will enable  us to prove Theorem \ref{thm:transversal-general}.
	
	\begin{theorem}\label{thm:general1}
		For  any  permutations $\alpha$, $\beta$ and $\tau$, if Peak is  equidistributed over  
		$\mathcal{T}_{\mu}(\alpha)$ and $\mathcal{T}_{\mu}(\beta)$ for any Young diagram $\mu$,
		then there exists a peak set preserving bijection $\Theta$ between $\mathcal{T}_{\lambda}(\alpha\oplus \tau)$ and $\mathcal{T}_{\lambda}(\beta\oplus \tau)$ for any Young diagram $\lambda$.
	\end{theorem}
	
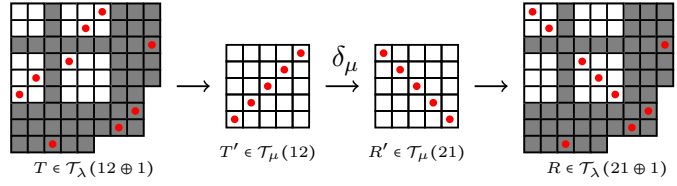
\begin{figure}[H]
		\begin{center}
			\begin{tikzpicture}[font =\small , scale = 0.22,line width = 0.7pt]
				\foreach \i / \j in {1/9,2/9,3/9,4/9,5/9,6/8,7/8,8/8,9/5}
				\foreach \k in {1,...,\j}
				{
					\filldraw[gray] (\i,-\k)rectangle(\i + 1,-\k -1);
					\draw (\i,-\k)rectangle(\i + 1,-\k -1);
				}
				\foreach \i / \j in {1/1,1/2,1/4,1/5,1/6,2/1,2/2,2/4,2/5,2/6,4/1,4/2,4/4,4/5,4/6,5/1,5/2,5/4,5/5,5/6,6/1,6/2,6/4,6/5,6/6}
				{   \filldraw[white] (\i,-\j)rectangle(\i + 1,-\j -1);
					\draw(\i,-\j)rectangle(\i + 1,-\j -1);
				}
				
				\foreach \i / \j in {1/6,2/5,3/9,4/4,5/2,6/1,7/8,8/7,9/3}
				\filldraw[red](\i+0.5,-\j-0.5)circle(5pt);
				\node[font = \tiny] at(6,-11){$T\in\mathcal{T}_{\lambda}(12\oplus 1)$};
				\draw[->](11,-6)--(13,-6);
				
				\tikzmath{\r = 13;\d = -2.5;};
				\foreach \i / \j in {1/5,2/5,3/5,4/5,5/5}
				\foreach \k in {1,...,\j}
				\draw (\i + \r,-\k + \d)rectangle(\i + 1 + \r,-\k -1 +\d);
				\foreach \i / \j in {1/5,2/4,3/3,4/2,5/1}
				\filldraw[red](\i+0.5 + \r,-\j-0.5 +\d)circle(5pt);
				\node[font = \tiny] at(16.5,-10){$T'\in\mathcal{T}_{\mu}(12)$};
				\draw[->](20,-6)--(22,-6)node[left=5,above]{$\delta_{\mu}$};
				
				\tikzmath{\r = 22;\d = -2.5;};
				\foreach \i / \j in {1/5,2/5,3/5,4/5,5/5}
				\foreach \k in {1,...,\j}
				\draw (\i + \r,-\k + \d)rectangle(\i + 1 + \r,-\k -1 +\d);
				\foreach \i / \j in {1/1,2/2,3/3,4/4,5/5}
				\filldraw[red](\i+0.5 + \r,-\j-0.5 +\d)circle(5pt);
				\node[font = \tiny] at(25.5,-10){$R'\in\mathcal{T}_{\mu}(21)$};
				\draw[->](29,-6)--(31,-6);
				
				\tikzmath{\r=31;};
				\foreach \i / \j in {1/9,2/9,3/9,4/9,5/9,6/8,7/8,8/8,9/5}
				\foreach \k in {1,...,\j}
				{   \filldraw[gray] (\i+\r,-\k)rectangle(\i + 1+\r,-\k -1);
					\draw (\i +\r,-\k)rectangle(\i + 1 + \r,-\k -1);
				}
				\foreach \i / \j in {1/1,1/2,1/4,1/5,1/6,2/1,2/2,2/4,2/5,2/6,4/1,4/2,4/4,4/5,4/6,5/1,5/2,5/4,5/5,5/6,6/1,6/2,6/4,6/5,6/6}
				{   \filldraw[white] (\i+\r,-\j)rectangle(\i + 1+\r,-\j -1);
					\draw(\i+\r,-\j)rectangle(\i + 1+\r,-\j -1);
				}
				
				\foreach \i / \j in {1/1,2/2,3/9,4/4,5/5,6/6,7/8,8/7,9/3}
				\filldraw[red](\i+0.5+ \r,-\j-0.5)circle(5pt);
				
				\node[font = \tiny] at(6+\r,-11){$R\in\mathcal{T}_{\lambda}(21\oplus 1)$};	
			\end{tikzpicture}
		\end{center}
		\caption{ An example of the bijection $\Theta$ between $\mathcal{T}_{\lambda}(\alpha\oplus \tau)$ and $\mathcal{T}_{\lambda}(\beta\oplus \tau)$.}\label{fig:Theta}
	\end{figure}

	\pf
	For any Young diagram $\mu$, let $\delta_{\mu}: \mathcal{T}_{\mu}(\alpha) \rightarrow  \mathcal{T}_{\mu}(\beta)$ be the peak set preserving bijection implied by the hypothesis. 
	Now we give a description of the bijection $\Theta: \mathcal{T}_{\lambda}(\alpha\oplus \tau) \rightarrow \mathcal{T}_{\lambda}(\beta\oplus \tau)$.
	Given a transversal $T=t_1t_2\cdots t_n \in \mathcal{T}_{\lambda}(\alpha\oplus \tau)$, color the
	square $(c,r)$ by white if the board of $\lambda$ lying below and to the right of it
	contains the pattern $\tau$, or gray otherwise.
	Find out all the $1$'s colored by gray and color the corresponding rows and columns by gray.
	Delete all the squares colored by gray, as well as the fillings in them.
	Let $\mu$ be the resulting board and let
	$T' = t_1't_2'\cdots t_k'$ be the corresponding $01$-filling of $\mu$.
	One can easily verify that $\mu$ is a Young diagram and
	$T'$ is a  transversal of shape $\mu$.
	Notice that $T$ avoids the pattern $\sigma\oplus \tau$ if and only if 
	$T'$ avoids the patter $\sigma$ for any permutation $\sigma$.
	Thus we have $T'$ avoids $\alpha$.
	By applying the bijection $\delta_{\mu}$, we  obtain a transversal
	$\delta_{\mu}(T')=R' =r_1'r_2'\cdots r_k' \in \mathcal{T}_{\mu}(\beta)$ with the same  peak set as $T'$. 
	Restoring the gray cells of $\lambda$ and  their contents,
	we obtain a  transversal $R = \Theta(T)$ of shape $\lambda$ which avoids the pattern $\beta\oplus \tau$.
 Figure \ref{fig:Theta} illustrates an example of $\Theta: \mathcal{T}_{\lambda}(\alpha\oplus \tau) \rightarrow \mathcal{T}_{\lambda}(\beta\oplus \tau)$ 
	where $\alpha = 12$, $\beta = 21$, $\tau = 1$ and $\lambda = (9,9,9,9,9,8,8,8,5)$.
	
	In order to show that the map   $\Theta$ is a bijection, it suffices to show that  the above procedure is invertible.  
	It is obvious that    $\Theta$ changes the $01$-filling located at the white squares  and leaves the $01$-filling located at the  gray  squares  fixed.   Hence  when  applying  the  inverse  map   $\Theta^{-1}$,   the  coloring  of $R$ will result in the same  Young diagram  $\mu$  and the same transversal $R'$  such that when  applying the inverse bijection   $\delta_{\mu}^{-1}$ to $R'$, we will recover the same transversal $T'$ and hence the same transversal $T$. 
	
	Now it is sufficient to show $\mathrm{Peak}(T) = \mathrm{Peak}(R)$.
	Assume that $i \in \mathrm{Peak}(T)$. 
	If the square $(i,t_i)$ is colored by gray, then by the coloring rule and the fact that $t_{i-1} < t_i > t_{i+1}$, the squares $(i-1,j)$ and $(i+1,j)$ with $j\geq t_i$ are also colored by gray.
	Then by the construction of $\Theta$, 
	one can easily check  that $i \in \mathrm{Peak}(R)$. 
	Now we assume that
	the square $(i,t_i)$ is colored by white and the lowest  white colored square in column $i$ is $(i,r)$. 
	Assume that column $i$ of $\lambda$ corresponds column
	 $i'$ in $\mu$.  
	Then the lowest white colored square in column $i-1$ (resp. $i+1$)
	is $(i-1,r)$ (resp. $(i+1,r)$) as $i\in \mathrm{Peak}(T)$. 
	This implies that $c_{i'-1}(\lambda')=c_{i'}(\lambda')=c_{i'+1}(\lambda')$ and the relative positions of $1$'s in columns $i'-1$, $i'$ and $i'+1$ in $T'$ are the same as those of  $1$'s in columns $i-1$, $i$ and $i+1$ in $T$.  Hence we have $i'\in  \mathrm{Peak}(T') =\mathrm{Peak}(R')$.
	By restoring the gray cells of $\lambda$ and  their contents, it is easily seen that $i$ is a peak of $R$.
	To conclude, we have that $\mathrm{Peak}(T) \subseteq \mathrm{Peak}(R)$.
	By similar arguments, one can show that 
	$\mathrm{Peak}(R) \subseteq \mathrm{Peak}(T)$.
	Hence we have concluded that $\mathrm{Peak}(T) = \mathrm{Peak}(R)$, completing the proof. 
	\qed
	
	Combining Theorems \ref{thm:transversal} and \ref{thm:general1},
	we are led to a proof of Theorem \ref{thm:transversal-general}.
	
	For symmetric transversals,  one can deduce the following analogue  of Theorem \ref{thm:general1} by making a slight adaption of the coloring rule in the proof. More precisely, we color  the square
	by white if the board  lying below and to the right of it contains the
	patterns $\tau$ or $\tau^{-1}$, or gray otherwise. 
	
	\begin{theorem} \label{thm:general2}
		Let $\alpha$ and $\beta$ be any involutions and let $\tau$ be any permutation.
		If the set-valued statistic  Peak is  equidistributed   over $\mathcal{ST}_{\mu}(\alpha)$
		and $\mathcal{ST}_{\mu}(\beta)$ for any self-conjugate 
		Young diagram $\mu$, 
		then there exists a peak set preserving bijection  between
		$\mathcal{ST}_{\lambda}(\alpha \oplus \tau)$ and $\mathcal{ST}_{\lambda}(\beta \oplus \tau)$ for any self-conjugate 
		Young diagram $\lambda$.
	\end{theorem}

	Combining Theorems \ref{thm:symm-transversal} and \ref{thm:general2}, we are led to a proof of 
	Theorem \ref{thm:symm-transversal-general}. 	To prove Theorem \ref{thm:AI}, we need the following lemma.
	
	\begin{lemma}\label{lem:AI}
		
		\begin{enumerate}[label=\upshape(\roman*)]
			\item For any permutations $\alpha$, $\beta$ and  any nonempty permutation $\tau$, if  the set-valued statistic Peak is equidistributed over  $\mathcal{T}_{\mu}(\alpha)$
			and $\mathcal{T}_{\mu}(\beta)$ for any 	Young diagram $\mu$,
			then we have  $|\mathcal{A}_n(\alpha\oplus\tau)|=|\mathcal{A}_n(\beta\oplus\tau)|$.	
			
			\item   Let  $\alpha$ and  $\beta$ be any involutions and let $\tau$ be   any nonempty permutation.   If the set-valued statistic Peak is equidistributed over   $\mathcal{ST}_{\mu}(\alpha)$
			and $\mathcal{ST}_{\mu}(\beta)$ for any self-conjugate 
			Young diagram $\mu$,
			then we have $|\mathcal{AI}_n(\alpha\oplus\tau)|=|\mathcal{AI}_n(\beta\oplus\tau)|$.
		\end{enumerate}
	\end{lemma}
	
	\pf
	We shall only prove (\upshape \rmnum{1}), as (\upshape \rmnum{2}) can be 
	deduced by similar arguments.
	By the hypothesis and Theorem \ref{thm:general1}, 
	the bijection $\Theta$  induces a peak set preserving between $\mathcal{S}_{n}(\alpha \oplus \tau)$ and 
	$\mathcal{S}_{n}(\beta \oplus \tau)$.
	Let $\pi\in \mathcal{A}_n(\alpha\oplus \tau) $
	and  $\Theta(\pi)=\sigma=\sigma_1\sigma_2\cdots \sigma_{n}$.
	Notice that the  permutation matrix of an alternating permutation $\pi \in \mathcal{A}_{n}$  is a transversal $T = t_1t_2\cdots t_n$ of a square Young diagram $\lambda=(\lambda_1, \lambda_2, \ldots, \lambda_n)$ with $\lambda_1=\lambda_2=\cdots=\lambda_n=n$  such that $\mathrm{Peak}(T)=\{2, 4, \ldots, 2\lfloor {n-1\over 2}\rfloor\}$ and when $n$ is even, $t_{n-1}<t_n$.
	In order to show that $\Theta(\pi)\in \mathcal{A}_n(\beta\oplus \tau) $,    
	it suffices to show that  $\sigma_{2k}>\sigma_{2k-1}$  when $n=2k$.
	This is justified by the fact that, when applying the bijection $\Theta$,  column $2k$  is always colored by gray and  all the squares colored by white (if any) in column $2k-1$  are positioned  above the $1$ located in column $2k$ as  $\tau$ is nonempty.
	Hence $\Theta$ induces a bijection between $\mathcal{A}_n(\alpha\oplus \tau)$ and $\mathcal{A}_n(\beta\oplus \tau)$.
	This completes the proof.
	\qed
	
	Combining  Theorem \ref{thm:transversal} and Lemma \ref{lem:AI}, 
	we derive that $|\mathcal{A}_n(I_k\oplus\tau)|=|\mathcal{A}_n(J_k\oplus\tau)|$
	for any nonempty permutation $\tau$ and for any positive integer $k$ which was first proved by 
	Yan \cite{Yan2013}.
	Similarly, by Theorem \ref{thm:symm-transversal} and Lemma \ref{lem:AI},
	we have $|\mathcal{AI}_n(I_k\oplus\tau)|=|\mathcal{AI}_n(J_k\oplus\tau)|$ for any nonempty pattern $\tau$ and for any positive integer $k$,
	completing the proof of Theorem \ref{thm:AI}.



	\section*{Acknowledgments}
	The work  was supported by
	the National Natural
	Science Foundation of China (12071440 and 11801378).
	
	\section*{Data availability statements}
	Data sharing is not applicable to this article as no datasets were generated or analyzed during the current study.


\begin{thebibliography}{100}
		
		\bibitem{Aguiar}
		M. Aguiar, K. Nyman, and R. Orellana, New results on the peak algebra, 
		{\em J. Algebraic Combin.}, {\bf 23} (2006), 149-–188.
		
		\bibitem{BW}
		E. Babson and J. West, The permutations $123p_4\ldots p_m$ and $321p_4\ldots p_m$ are Wilf-equivalent, {\em Graphs Combin.}, {\bf 16} (2000),
		373--380.
		
		\bibitem{BWX}
		J. Backelin, J. West, and G. Xin.  \newblock Wilf-equivalence for singleton classes.  {\em Adv. Appl.
			Math.},   {\bf 38} (2007), 133--148.
		
		\bibitem{Barnabei2011}
		M. Barnabei, F. Bonetti, and M. Silimbani, Restricted involutions
		and Motzkin paths, {\em Adv.  Appl. Math.}, {\bf 47} (2011),
		102--115.
		
		\bibitem{Barnabei}
		M. Barnabei, F. Bonetti, N. Castronuovo and M. Silimbani, Pattern avoiding alternating involutions,
		{\em ECA}, {\bf 3:1} (2023), Article \#S2R4.
		
		\bibitem{Billera}
		L.J. Billera, S.K. Hsiao, and S. van Willigenburg, Peak quasisymmetric functions and Eulerian enumeration, {\em Adv. Math.}, {\bf 176} (2003), 248–-276.
		
		
		
		\bibitem{Bloom1}
		J. Bloom, D. Saracino,   On bijections for pattern-avoiding permutations, {\em J. Combin. Theory Ser. A},   {\bf 116}(2009), 1271--1284. 
		\bibitem{Bloom2}
		J. Bloom, D. Saracino,  Another look at bijections for pattern-avoiding permutations,    {\em Adv. Appl. Math.},   {\bf 45}(2010), 395--409,   2010.
		
		\bibitem{Bloom3}
		J. Bloom,   A refinement of Wilf-equivalence for patterns of length 4,    {\em J. Combin. Theory Ser. A},    {\bf 124} (2014), 166--177.
		
		
		
		\bibitem{Bona}
		M. B{\'o}na, On a family of conjectures of Joel Lewis on alternating permutations, {\em Graphs Combin.}, {\bf 30} (2014), 521--526.
		
		\bibitem{Bona2016}
		M. B\'ona, C. Homberger, J. Pantone, and V. Vatter,  Pattern-avoiding
		involutions: exact and asymptotic enumeration, {\em Australas. J. Combin.}, {\bf 64} (2016), 88--119.
		
		
		
		\bibitem{Bousquet}
		M. Bousquet-M\'elou and  E. Steingr\'imsson, Decreasing subsequences in permutations and Wilf equivalence
		for involutions, {\em J. Algebraic Combin.}, {\bf 22} (2005), 383--409.
		
		\bibitem{Chan}
		J.H.C. Chan, 
		An infinite family of inv-Wilf-equivalent
		permutation pairs, {\em European J. Combin.}, {\bf 44} (2015), 57--76.
		
		\bibitem{ChenD}
		W.Y.C. Chen, E.Y.P. Deng, R.R.X. Du, R.P. Stanley, and C.H. Yan, Crossings and nestings of matchings and partitions, {\em Trans. Amer. Math. Soc.}, {\bf 359} (2007), 1555--1575.
		
		
		\bibitem{Chen}
		J.N. Chen, W.Y.C. Chen, and R.D.P. Zhou, On pattern avoiding alternating permutations, {\em European J. Combin.}, {\bf 40} (2014), 11--25.
		
		
		
		\bibitem{dd}
		T.  Dokos,  T.  Dwyer, B. P. Johnson,  B.E. Sagan, and  K. Selsor,   Permutation patterns and statistics, {\em
			Discrete Math.},    {\bf 312}(2012),  2760--2775. 
		
		\bibitem{Dukes}
		W.M.B. Dukes, V. Jel\'inek, T. Mansour, and A. Reifegerste, New equivalences for pattern avoiding involutions,
		{\em Proc. Amer. Math. Soc.}, {\bf 137} (2009), 457--465
		
		
		\bibitem{Guibert}
		O. Guibert, Combinatoire des permutations \` a motifs exclus, en liaison avec mots, cartes planaires et tableaux de Young,  PhD thesis, LaBRI, Universit\'e Bordeaux 1, 1995.
		
		
		
		
		\bibitem{Jaggard}
		A.D. Jaggard, Prefix exchanging and pattern avoidance by involutions, {\em Electron.  J. Combin.}, {\bf 9(2)} (2003), R16.
		
		
		
		
		\bibitem{Knuth}
		D.E. Knuth, The art of computer programming, {\em sorting and searching}, vol. 3, Addison-Wesley, 1973.
		
		
		\bibitem{Lewis2009}
		J.B. Lewis, Alternating, pattern-avoiding permutations, {\em Electron. J. Combin.}, {\bf 16} (2009), N7.
		
		\bibitem{Lewis2011}
		J.B. Lewis, Pattern avoidance for alternating permutations and Young tableaux, {\em J. Combin. Theory Ser. A}, {\bf 118} (2011), 1436--1450.
		
		\bibitem{Lewis2012}
		J.B. Lewis, Generating trees and pattern avoidance in alternating permutations, {\em Electron. J. Combin.}, {\bf 19} (2012), P21.
		
		\bibitem{Mansour2003}
		T. Mansour, Restricted 132-alternating permutations and Chebyshev polynomials, {\em Ann.  Combin.}, {\bf 7} (2003), 201--227.
		
		\bibitem{Petersen}
		T.K. Petersen, 	Enriched P-partitions and peak algebras,
		{\em Adv. Math.}, {\bf 209} (2007), 61–-610.
		
		\bibitem{Schocker}
		M. Schocker, The peak algebra of the symmetric group revisited,
		{\em Adv. Math.}, {\bf 192}	(2005), 259-–309.
		
		
		\bibitem{Simion}
		R. Simion and F. Schmidt, Restricted permutations, {\em European J. Combin.}, {\bf 6} (1985), 383--406.
		
		
		\bibitem{Stanley}
		R.P. Stanley, Alternating permutations and symmetric functions,
		{\em J. Combin. Theory Ser. A}, {\bf 114} (2007), 436--460.
		
		\bibitem{StanleyVol2}
		R.P. Stanley, Enumerative combinatorics, vol. 2, Cambridge University Press, Cambridge, 1999.
		
		\bibitem{Stembridge}
		J.R. Stembridge, Enriched P-partitions, {\em Trans. Amer. Math. Soc.}, {\bf 349} (1997), 763-–788.
		
		\bibitem{Strehl}
		V. Strehl, Enumeration of alternating permutations according to peak sets, {\em J. Combin. Theory Ser. A}, {\bf 24} (1978), 238-–240.
		
		\bibitem{Sundaram}
		S. Sundaram, The Cauchy identity for $Sp(2n)$,   {\em J. Combin. Theory Ser. A}, {\bf 53} (1990), 209--238.
		
		\bibitem{Warren}
		D. Warren and E. Seneta, Peaks and Eulerian numbers in a random sequence, {\em J. Appl. Probab.},  {\bf 33} (1996), 101-–114.
		
		
		\bibitem{xinzhang}
		G. Xin and T.Y.J. Zhang, Enumeration of bilaterally symmetric 3-noncrossing partitions, {\em Discrete Math.},  {\bf 309} (2009), 2497--2509.
		
		
		\bibitem{Yan2012}
		Y.X. Xu and S.H.F. Yan, Alternating permutations with restrictions and standard Young tableaux, {\em Electron. J. Combin.}, {\bf 19(2)} (2012), P49.
		
		
		\bibitem{Yan2013}
		S.H.F. Yan, On  Wilf equivalence for alternating permutations, {\em Electron. J. Combin.}, {\bf 20(3)} (2013), P58.
		
		\bibitem{Yan2015}
		S.H.F. Yan,  H. Ge, and Y. Zhang, On a refinement of Wilf-equivalence for permutations, {\em Electron. J. Combin.}, {\bf 22(1)} (2015),  P1.20. 
		
		\bibitem{Yan2023}
		S.H.F. Yan, L.T. Wang, and R.D.P. Zhou, On refinements of Wilf-equivalence for involutions,
		{\em J. Algebraic Combin.}, accepted, arXiv:2212.01800v1 [math.CO]. 
		
		\bibitem{Zhou}
		R.D.P. Zhou, S.H.F. Yan,  Combinatorics of exterior peaks on pattern avoiding  symmetric transversals, {\em preprint}. 
		
	\end{thebibliography}
\end{document}